\documentclass[12pt,reqno]{amsart}
\usepackage{amsmath}
\usepackage{amssymb}
\usepackage{amsfonts}
\usepackage{eucal}
\usepackage{epsfig}

\usepackage{waveletnotation}

\newcommand{\B}{{\mathfrak{B}}}
\newcommand{\seqspace}{{\mathfrak{b}}}
\newcommand{\Wop}{\mathcal{W}}

\newtheorem{lemma}{Lemma}
\newtheorem{prop}[lemma]{Proposition}
\newtheorem{cor}[lemma]{Corollary}
\newtheorem{theorem}[lemma]{Theorem}

\begin{document}

\title[Frequency-based Nonhomogeneous Wavelet Systems in the Distribution Space]{Pairs of Frequency-based Nonhomogeneous Dual Wavelet Frames in the Distribution Space}

\author{Bin Han}

\thanks{Research supported in part by NSERC Canada under Grant
RGP 228051. \hfill  December 18, 2009.
}

\address{Department of Mathematical and Statistical Sciences,
University of Alberta, Edmonton,\quad Alberta, Canada T6G 2G1.
\quad {\tt bhan@math.ualberta.ca}\quad
{\tt http://www.ualberta.ca/$\sim$bhan}
}

\makeatletter \@addtoreset{equation}{section} \makeatother
\begin{abstract}
In this paper, we study nonhomogeneous wavelet systems which have close relations to the fast wavelet transform and homogeneous wavelet systems. We introduce and characterize a pair of frequency-based nonhomogeneous dual wavelet frames in the distribution space; the proposed notion enables us to completely separate the perfect reconstruction property of a wavelet system from its stability property in
function spaces. The results in this paper lead to a natural explanation for the oblique extension principle, which has been widely used to construct dual wavelet frames from refinable functions, without any a priori condition on the generating wavelet functions and refinable functions.
A nonhomogeneous wavelet system,
which is not necessarily derived from refinable functions via a multiresolution analysis, not only has a natural multiresolution-like structure that is closely linked to the fast wavelet transform, but also plays a basic role in understanding many aspects of wavelet theory.
To illustrate the flexibility and generality of the approach in this paper, we further extend our results to nonstationary wavelets with real dilation factors and to nonstationary wavelet filter banks having
the perfect reconstruction property.
\end{abstract}


 \maketitle

\bigskip

\pagenumbering{arabic}


\section{Introduction and Motivations}

In wavelet analysis, we often use translation, dilation, and modulation of functions. For a function $f: \R \rightarrow \C$, throughout the paper we shall use the following notation
\begin{equation}\label{tdm}
f_{\gl; k, n}(x):=|\gl|^{1/2} e^{-\iu n \gl x} f(\gl x-k) \quad \mbox{and}\quad
f_{\gl; k}:=f_{\gl;k,0}=|\gl|^{1/2} f(\gl \cdot
-k),
 \qquad x, \gl, k, n\in \R,
\end{equation}
where $\iu$ denotes the imaginary unit. In this paper we shall use $\df \in \R \bs \{0\}$ as a dilation factor. In applications, $\df$ is often taken to be a positive integer greater than one, in particular, the simplest case $\df=2$ is often used.

Classical wavelets are often defined and studied in the time/space domain with the generating wavelet functions belonging to the square integrable function space $\Lp{2}$.
For $\df\in \R\bs \{0\}$ and for a subset $\Psi$ of square integrable functions in $\Lp{2}$, linked to discretization of a continuous wavelet transform (see \cite{Chui:book,Daub:book,Mallat:book,Meyer:book}), the following homogeneous wavelet system
\begin{equation}\label{WS:hom}
\WS(\Psi):=
\{\psi_{\df^j; k}\; : \; j\in \Z, k\in \Z, \psi\in \Psi\}
\end{equation}
is generated by the translation and dilation of the wavelet functions in $\Psi$ and has been extensively studied in the function space $\Lp{2}$ in the literature of wavelet analysis. To mention only a few references here, see \cite{BGN}--\cite{RonShen:twf}. In this paper, however, we shall see that it is more natural to study a nonhomogeneous wavelet system in the frequency domain. It is important to point out here that the elements in a set $S$ of this paper are not necessarily distinct and $h \in S$ in a summation means that $h$ visits every element (with multiplicity) in $S$ once and only once.
For example, for $\Psi=\{\psi^1, \ldots, \psi^\mpsi\}$, all the functions $\psi^1, \ldots, \psi^\mpsi$ are not necessarily distinct and $\psi\in \Psi$ in \eqref{WS:hom} means $\psi=\psi^1, \ldots, \psi^\mpsi$.

Most known classical homogeneous wavelet systems $\WS(\Psi)$ in the literature are often derived from scalar refinable functions $\phi$ or from refinable function vectors $\vec{\phi}=(\phi^1, \ldots, \phi^\mphi)^T$ in $\Lp{2}$ (\cite{Chui:book,Daub:book,Mallat:book,Meyer:book}). Let us recall the definition of a scalar refinable function here. A function or distribution $\phi$ on $\R$ is said to be {\it refinable} (or $\df$-refinable) if there exists a sequence $a=\{a(k)\}_{k\in \Z}$ of complex numbers, called the {\it refinement mask} or the {\it low-pass filter} for the scalar refinable function $\phi$, such that
\begin{equation}\label{refeq}
\phi=|\df|\sum_{k \in \Z} a(k) \phi(\df \cdot-k)
\end{equation}
with the above series converging in a proper sense, e.g., in $\Lp{2}$. Wavelet functions in the generating set $\Psi$ of a homogeneous wavelet system $\WS(\Psi)$ are often derived from the refinable function $\phi$ by
\begin{equation}\label{wavelet}
\psi=|\df|\sum_{k \in \Z} b^{\psi}(k) \phi(\df \cdot-k), \qquad \psi\in \Psi,
\end{equation}
where $b^\psi=\{b^\psi(k)\}_{k\in \Z}$ are sequences on $\Z$, called {\it wavelet masks} or {\it high-pass filters}. For the infinite series in \eqref{refeq} and \eqref{wavelet} to make sense, one often imposes some decay condition on the refinable function $\phi$ and wavelet filters $a, b^\psi$ so that all the infinite series in \eqref{refeq} and \eqref{wavelet} are well-defined in a proper sense. Nevertheless, even for the simplest case of a compactly supported scalar refinable function (or distribution) $\phi$ with a finitely supported mask $a$, the associated refinable function $\phi$ with mask $a$ does not always belong to $\Lp{2}$. In fact, it is far from trivial to check whether $\phi\in \Lp{2}$ in terms of its mask $a$, see
\cite{Han:sm,Han:acm:2006} and references therein for detail.
One of the motivations of this paper is to study wavelets and framelets without such stringent conditions on either the generating wavelet functions $\phi, \psi$ or their wavelet filters $a, b^\psi$ for $\psi\in \Psi$.

For $f\in \Lp{1}$, the Fourier transform used in this paper is defined to be $\hat f(\xi):=\int_\R f(x) e^{-ix\xi} dx$, $\xi\in \R$ and can be naturally extended to square integrable functions and tempered distributions.
Under certain assumptions, taking Fourier transform on both sides of \eqref{refeq} and \eqref{wavelet}, one can easily rewrite  \eqref{refeq} and \eqref{wavelet} in the frequency domain as follows:
\begin{equation}\label{refeq:F}
\hat \phi(\df \xi)=\hat a(\xi) \hat \phi(\xi), \qquad a.e.\; \xi\in \R
\end{equation}
and
\begin{equation}\label{wavelet:F}
\hat{\psi}(\df \xi)=\wh{b^{\psi}}(\xi) \hat \phi(\xi), \qquad a.e.\; \xi\in \R, \; \psi\in \Psi
\end{equation}
provided that all the $2\pi$-periodic (Lebesgue) measurable functions $\hat a(\xi)=\sum_{k\in \Z} a(k) e^{-\iu k\xi}$ and similarly $\wh{b^\psi}$ are properly defined. In the following, we shall see that it is often more convenient to work with \eqref{refeq:F} and \eqref{wavelet:F} in the frequency domain rather than \eqref{refeq} and \eqref{wavelet} in the time/space domain.
If there exist positive real numbers $\tau$ and $C$ such that the $2\pi$-periodic measurable function $\hat a$ satisfies
$|1-\hat a(\xi)| \le C |\xi|^\tau$ for almost every $\xi\in [-\pi, \pi]$ (this condition is automatically satisfied with $\tau=1$ if $\hat a$ is a $2\pi$-periodic trigonometric polynomial with $\hat a(0)=1$), for a dilation factor $\df$ such that $|\df|>1$, then it is easy to see (also c.f. section~3) that one can define a measurable function $\fphi$ such that
\begin{equation}\label{refeq:standard}
\fphi(\xi):=\prod_{j=1}^\infty \hat a(\df^{-j}\xi), \qquad a.e. \; \xi\in \R.
\end{equation}
%
Regardless of whether $\fphi$ in \eqref{refeq:standard} is a square integrable function or not, $\fphi$ is a well-defined measurable function obviously satisfying the frequency-based refinement equation \eqref{refeq:F} with $\hat \phi$ being replaced by $\fphi$. The function $\fphi$ is called the (frequency-based) {\it standard refinable function} with mask $\hat a$ and dilation $\df$.
All the wavelet functions $\hat{\psi}$ in \eqref{wavelet:F} with $\hat \phi=\fphi$ are also well-defined measurable functions provided that all $\wh{b^{\psi}}$ are measurable. This motivates us to study refinable functions and wavelets in the frequency domain using \eqref{refeq:F} and \eqref{wavelet:F} so that we can avoid some technical issues such as the convergence of the infinite series in \eqref{refeq} and \eqref{wavelet} as well as membership in $\Lp{2}$ of the generating refinable function $\phi$ and the generating wavelet functions $\psi \in \Psi$.

For wavelets derived from refinable functions or refinable function vectors,
one of the most important key features of wavelets and framelets is its associated fast wavelet transform, which is based on the following nonhomogeneous wavelet system:
\begin{equation}\label{WS}
\WS_J (\Phi; \Psi):=\{ \phi_{\df^J;k}\; : \; k\in \Z, \phi\in \Phi\}\cup\{\psi_{\df^j;k}\; : \; j\ge J, k\in \Z, \psi\in \Psi\},
\end{equation}
where 
$J$ is an integer, representing the coarsest decomposition (or scale) level of its fast wavelet transform. In fact, a one-level fast wavelet transform is just a transform between two sets of wavelet coefficients of a given function represented under two nonhomogeneous wavelet systems at two consecutive scale levels. Naturally, for a multi-level wavelet transform, there is a underlying sequence of nonhomogeneous wavelet systems at all scale levels, instead of just one single wavelet system. For a given 
$\WS_{J_0} (\Phi; \Psi)$ at some scale $J_0$, we shall see in this paper that via the dilation operation it naturally produces a sequence of nonhomogeneous wavelet systems
$\WS_J (\Phi; \Psi)$ for all integers $J$ with almost all properties preserved. Consequently, it often suffices to study only one nonhomogeneous wavelet system instead of a sequence of them.  This desirable property of nonhomogeneous wavelet systems is not shared by homogeneous wavelet systems. Furthermore, as $J\to -\infty$, the limit of the sequence $\{\WS_J (\Phi; \Psi)\}_{J\in \Z}$ will naturally lead to a homogeneous wavelet system $\WS(\Psi)$. Hence, in certain sense, a homogeneous wavelet system could be regarded as the limit system of a
sequence of nonhomogeneous wavelet systems. See section~3 for more detail.

For a homogeneous wavelet system $\WS(\Psi)$ that is derived from a refinable function or a refinable function vector, due to the absence of a refinable function $\phi$ in the system, the homogeneous wavelet system $\WS(\Psi)$ does not automatically correspond to a fast wavelet transform without ambiguity. In fact, the wavelet functions in $\Psi$ of $\WS(\Psi)$ could be derived from many other (equivalent) refinable functions, which correspond to different fast wavelet transforms with different sets of wavelet filters. More precisely, for a $2\pi$-periodic measurable function $\theta$ such that $\theta(\xi)\ne 0$ for almost every $\xi\in \R$, define
$\hat{\mathring{\phi}}(\xi):=\theta(\xi)\hat \phi(\xi)$, then it is evident that $\hat{\mathring{\phi}}$ is also refinable and satisfies
\[
\hat{\mathring{\phi}}(\df\xi)=[\theta(\df\xi)\hat a(\xi)/\theta(\xi)] \hat {\mathring{\phi}}(\xi) \quad \hbox{and}\quad \hat{\psi}(\df\xi)=[\wh{b^\psi}(\xi)/\theta(\xi)] \hat{\mathring{\phi}}(\xi), \qquad a.e.\; \xi\in \R, \; \psi\in \Psi.
\]
Such a change of generators from a refinable function $\phi$ to another equivalent refinable function $\mathring{\phi}$ is in fact the key idea in the oblique extension principle (OEP) in \cite{CHS,DH:dwf,DHRS,Han:oep} (also see
\cite{E1,E2,Han:jcam,HanMo:siamaa,HanMo:acha,HanShen:siam:2008,HanShen:ca:2009}) to construct compactly supported homogeneous wavelet frames $\WS(\Psi)$ in $\Lp{2}$ with high vanishing moments derived from refinable functions and refinable function vectors. See \cite{DHRS,Han:oep} for a detailed discussion on a fast wavelet transform based on a homogeneous wavelet system obtained via OEP from refinable function vectors. The effect of the change of generators on its fast wavelet transform is addressed in \cite{Han:oep}.
As we shall see in sections~3 and 4,
nonhomogeneous wavelet systems are closely related to nonstationary wavelets (see \cite{CHS:ntwf,CohenDyn,HanShen:siam:2008}) and are naturally employed in a pair of nonhomogeneous dual wavelet frames in a pair of dual Sobolev spaces introduced in \cite{HanShen:ca:2009}.
Due to these and other considerations, it seems more natural and more important for us to study nonhomogeneous wavelet systems $\WS_J (\Phi; \Psi)$ rather than the  extensively studied homogeneous wavelet systems $\WS(\Psi)$.
This allows us to understand better many aspects of wavelet theory such as a wavelet filter bank induced by OEP and its associated wavelets in the function setting without a priori condition on the generating wavelet functions.

Following the standard notation, we denote by $\D$ the linear space of all compactly supported $C^\infty$ (test) functions with the usual topology, and $\Dpr$ denotes the linear space of all distributions, that is, $\Dpr$ is the dual space of $\D$.
By duality, the definition in \eqref{tdm} for translation, dilation and modulation can be easily generalized from functions to distributions. Moreover, when $f$ in \eqref{tdm} is a square integrable function or more generally a tempered distribution, for $\gl\ne 0$, we have
\begin{equation}\label{F:tdm}
\wh{f_{\gl; k,n}}=e^{-\iu kn} \hat f_{\gl^{-1}; -n, k} \quad \mbox{and}\quad \wh{f_{\gl; k}}=\hat f_{\gl^{-1}; 0, k}.
\end{equation}
In this paper, we shall use boldface letters to denote functions/distributions (e.g., $\ff,\fg, \fphi, \fpsi, \fa, \fb$) or sets of functions/distributions (e.g., $\fPhi, \fPsi$) in the frequency domain.

Let $\fPhi$ and $\fPsi$ be two sets of distributions on $\R$.
For an integer $J$ and $\df\in \R \bs \{0\}$, we define {\it a frequency-based nonhomogeneous wavelet system $\FWS_J(\fPhi; \fPsi)$} to be
\begin{equation}\label{FWS}
\FWS_J(\fPhi; \fPsi):=\{\fphi_{\df^{-J};0, k}\; : \; k\in \Z, \fphi\in \fPhi\}\cup\{\fpsi_{\df^{-j}; 0, k}\; : \; j\ge J, k\in \Z, \fpsi\in \fPsi\}.
\end{equation}
Similarly, {\it a frequency-based homogeneous wavelet system} is defined as follows:
\begin{equation}\label{FHWS}
\FWS(\fPsi):=\{\fpsi_{\df^j;0,k}\, : \, j\in \Z, k\in \Z, \fpsi\in \fPsi\}.
\end{equation}
By \eqref{F:tdm}, it is straightforward to see that under the Fourier transform, the images
of $\WS_J (\Phi; \Psi)$ and $\WS(\Psi)$ with $\Phi, \Psi\subseteq \Lp{2}$ are simply $\FWS_J(\fPhi; \fPsi)$ and $\FWS(\fPsi)$, respectively, where
$\fPhi:=\{ \hat \phi \; : \; \phi\in \Phi\}$ and
$\fPsi:=\{ \hat \psi \; : \; \psi\in \Psi\}$.

For $1\le p \le \infty$, by $\lLp{p}$ we denote the linear space of all measurable functions $f$ such that $\int_K |f|^p<\infty$ for every compact subset $K$ of $\R$, with the usual modification for $p=\infty$ saying that $f$ is essentially bounded over $K$. Note that $\lLp{1}$ is just the set of all measurable functions that can be globally identified as distributions. So, $\lLp{1}$ is the most natural space for us to study wavelets and framelets in the distribution space.
For $f\in \Lp{2}$, it is evident that $\hat f\in \Lp{2}\subseteq \lLp{2}\subseteq \lLp{1}\subseteq \Dpr$. However, a distribution $\fpsi$ may not be a tempered distribution and therefore, Fourier transform may not be applied so that $\fpsi=\hat \psi$ holds for some tempered distribution $\psi$. Under the setting of tempered distributions on which the Fourier transform can apply, although all the definitions and results of this paper in the frequency domain could be equivalently translated into the time/space domain by the inverse Fourier transform on tempered distributions, to avoid notational complexity and to avoid the a priori underlying assumption that $f$ is a tempered distribution if the notion $\hat f$ is used, it seems very natural and convenient for us to work in the frequency domain in this paper.

For $\ff\in \D$ and $\fpsi\in \lLp{1}$, we shall use the following paring
\begin{equation}\label{pair}
\la \ff, \fpsi\ra:=\int_{\R} \ff(\xi) \overline{\fpsi(\xi)} d\xi \quad \mbox{and}\quad
\la \fpsi, \ff\ra:=\ol{\la \ff, \fpsi\ra}=
\int_{\R} \fpsi(\xi) \overline{\ff(\xi)} d\xi.
\end{equation}
When $\ff\in \D$ and $\fpsi\in \Dpr$, the duality pairings $\la \ff, \fpsi\ra$ and $\la \fpsi, \ff\ra$ are understood similarly as $\la \ff, \fpsi\ra:=\overline{\la \fpsi, \ff\ra}:=\overline{\fpsi(\overline{\ff})}$.
Now we are ready to introduce the key notion in this paper. Let
\begin{equation}\label{gen:set}
\fPhi=\{\fphi^1, \ldots, \fphi^\mphi\}, \quad \fPsi=\{\fpsi^1, \ldots, \fpsi^\mpsi\} \quad \mbox{and}\quad
\tilde \fPhi=\{\tilde \fphi^1, \ldots, \tilde\fphi^\mphi\}, \quad \tilde\fPsi=\{\tilde\fpsi^1, \ldots, \tilde\fpsi^\mpsi\}
\end{equation}
be subsets of $\Dpr$. Let $J\in \Z$ and $\df \in \R \bs \{0\}$,
we say that the pair $(\FWS_J(\fPhi; \fPsi), \FWS_J(\tilde \fPhi; \tilde \fPsi))$, where $\FWS_J(\fPhi; \fPsi)$ is defined in \eqref{FWS},
%
%
forms {\it a pair of frequency-based nonhomogeneous dual wavelet frames in the distribution space $\Dpr$} if the following identity holds
\begin{equation}\label{dwf}
\begin{split}
\sum_{\ell=1}^\mphi\sum_{k\in \Z} \la \ff, \fphi^\ell_{\df^{-J};0, k}\ra \la \tilde \fphi^\ell_{\df^{-J};0,k}, \fg\ra&+
 \sum_{j=J}^\infty  \sum_{\ell=1}^{\mpsi} \sum_{k\in \Z}
\la \ff, \fpsi^\ell_{\df^{-j};0,k}\ra \la \tilde \fpsi^{\ell}_{\df^{-j};0,k}, \fg\ra\\
&= 2\pi \la \ff, \fg\ra \qquad \forall\; \ff, \fg\in \D,
\end{split}
\end{equation}
where the infinite series in \eqref{dwf} converge in the following sense:
\begin{enumerate}
\item For every $\ff, \fg\in \D$, the following series
\begin{equation}\label{series}
\sum_{k\in \Z} \la \ff, \fphi^\ell_{\df^{-J};0,k}\ra \la \tilde \fphi^{\ell}_{\df^{-J};0,k}, \fg\ra\quad \mbox{and}\quad \sum_{k\in \Z} \la \ff, \fpsi^{\ell'}_{\df^{-j};0,k}\ra \la \tilde \fpsi ^{\ell'}_{\df^{-j};0,k}, \fg\ra
\end{equation}
converge absolutely for all integers $j\ge J$, $\ell=1, \ldots, \mphi$, and $\ell'=1, \ldots, \mpsi$.

\item For every $\ff, \fg\in \D$, the following limit exists and
\begin{equation}\label{def:dwf}
\lim_{J'\to +\infty} \left(\sum_{\ell=1}^{\mphi} \sum_{k\in \Z} \la \ff, \fphi^{\ell}_{\df^{-J};0,k}\ra \la \tilde \fphi^{\ell}_{\df^{-J};0,k}, \fg\ra+
 \sum_{j=J}^{J'-1} \sum_{\ell=1}^{\mpsi}\sum_{k\in \Z}
\la \ff, \fpsi^{\ell}_{\df^{-j};0,k}\ra \la \tilde \fpsi^{\ell}_{\df^{-j};0,k}, \fg\ra\right)=2\pi \la \ff, \fg\ra.
\end{equation}
\end{enumerate}
As we shall discuss in section~3, the above introduced notion
enables us to completely separate the perfect reconstruction property in \eqref{dwf} from its stability property in function spaces. Since the test function space $\D$ is dense in many function spaces, one could extend the perfect reconstruction property (or ``wavelet expansion'') in \eqref{dwf} to other function spaces, provided that the involved wavelet systems have stability in these function spaces.
Let us give a simple example here to illustrate this connection for the particular function space $\Lp{2}$.
Let $\fPhi$ and $\fPsi$ in \eqref{gen:set} be two subsets of distributions in $\Dpr$. We shall see in section~3 that $(\FWS_J(\fPhi; \fPsi), \FWS_J(\fPhi; \fPsi))$ forms a pair of frequency-based nonhomogeneous dual wavelet frames in the distribution space, if and only if,
$\fPhi, \fPsi$ are subsets of $\Lp{2}$ and
\begin{equation}\label{F:twf}
\sum_{\ell=1}^{\mphi} \sum_{k\in \Z} |\la \ff, \fphi^{\ell}_{\df^{-J};0,k}\ra|^2 + \sum_{j=J}^\infty  \sum_{\ell=1}^{\mpsi} \sum_{k\in \Z}
|\la \ff, \fpsi^{\ell}_{\df^{-j};0,k}\ra|^2= 2\pi \| \ff\|^2_{\Lp{2}} \qquad \forall\; \ff\in \Lp{2}.
\end{equation}
Moreover, when $|\df|>1$, as a direct consequence of \eqref{F:twf}, one automatically has
\[
\sum_{j\in \Z}  \sum_{\ell=1}^{\mpsi} \sum_{k\in \Z}
|\la \ff, \fpsi^{\ell}_{\df^{j};0,k}\ra|^2= 2\pi \|\ff \|^2_{\Lp{2}} \qquad \forall\; \ff \in \Lp{2}.
\]
Nonhomogeneous wavelet systems also have a close relation to refinable functions and refinable function vectors. Suppose that $\frac{1}{\sqrt{2\pi}}\FWS_J(\fPhi; \fPsi)$ (that is, multiply every element in $\FWS_J(\fPhi; \fPsi)$ by the factor $\frac{1}{\sqrt{2\pi}}$) is an orthonormal basis of $\Lp{2}$.
Denote $\vec{\fphi}:=(\fphi^1, \ldots, \fphi^{\mphi})^T$.
Without assuming in advance that $\vec{\fphi}$ is a refinable function vector and all $\fpsi^1, \ldots, \fpsi^{\mpsi}$ are derived from $\vec{\fphi}$, we can deduce that
$\vec{\fphi}$ must be a refinable function vector and all $\fpsi^1, \ldots, \fpsi^{\mpsi}$ must be derived from $\vec{\fphi}$ via similar relations as in \eqref{refeq:F} and \eqref{wavelet:F}.
See section~3 for more detail.

To have some rough ideas about our results on nonhomogeneous wavelet systems, here we present two typical results.
The following result is a special case of Theorem~\ref{thm:main:1:general}.

\begin{theorem}\label{thm:main:1}
Let $\df$ be an integer such that $|\df|>1$.
Let $\fPhi, \fPsi, \tilde \fPhi, \tilde \fPsi$ in \eqref{gen:set} be subsets of $\lLp{2}$. Then $(\FWS_J(\fPhi; \fPsi), \FWS_J(\tilde \fPhi; \tilde \fPsi))$ forms a pair of frequency-based nonhomogeneous dual wavelet frames in the distribution space $\Dpr$ for some integer $J$ (or for all integers $J$), if and only if, the following three statements hold
\begin{enumerate}
\item[{\rm (i)}] For all integers $k\in \Z$,
\begin{equation}\label{main:eq2}
\sum_{\ell=1}^\mphi \ol{\fphi^\ell(\df \xi)} \tilde \fphi^\ell(\df (\xi+2\pi k))+\sum_{\ell=1}^\mpsi \ol{\fpsi^\ell(\df\xi)}\tilde \fpsi^\ell(\df(\xi+2\pi k)) \\
=\sum_{\ell=1}^\mphi \ol{\fphi^\ell(\xi)}\tilde \fphi^\ell(\xi+2\pi k), \; a.e.\; \xi\in \R;
\end{equation}
\item[{\rm (ii)}] For all integers $k_0\in \Z \bs [\df \Z]$,
\begin{equation}\label{main:eq3}
\sum_{\ell=1}^\mphi \ol{\fphi^\ell(\xi)}\tilde \fphi^\ell(\xi+2\pi k_0)+\sum_{\ell=1}^\mpsi
\ol{\fpsi^\ell(\xi)}\tilde \fpsi^\ell(\xi+2\pi k_0)=0, \qquad a.e.\; \xi\in \R;
\end{equation}
\item[{\rm (iii)}] The following identity holds in the sense of distributions:
\begin{equation}\label{main:eq1}
\lim_{j\to +\infty} \sum_{\ell=1}^{\mphi} \ol{\fphi^\ell(\df^{-j}\cdot)}\tilde \fphi^\ell (\df^{-j}\cdot)=1,
\end{equation}
more precisely, $\lim_{j\to +\infty} \sum_{\ell=1}^{\mphi} \la \ol{\fphi^\ell(\df^{-j}\cdot)}\tilde \fphi^\ell (\df^{-j}\cdot), \ff\ra=\la 1, \ff\ra$ for all $\ff\in \D$.
\end{enumerate}
\end{theorem}

In the following, we make some remarks about Theorem~\ref{thm:main:1}.
We assumed in Theorem~\ref{thm:main:1} that all the generating functions in $\fPhi, \fPsi, \tilde \fPhi, \tilde \fPsi$ are from the space $\lLp{2}$. Note that $\lLp{2}$ includes the Fourier transforms of all compactly supported distributions and of all elements in all Sobolev spaces. This assumption on membership in $\lLp{2}$ can be weakened and is only used to guarantee the absolute convergence of the infinite series in \eqref{series}. See the remark after Lemma~\ref{lem:converg} in section~2 for more detail on this natural assumption.

If we assume additionally that $\ol{\fphi^\ell}\tilde \fphi^\ell \in \Lp{\infty}$ for all $\ell=1, \ldots, \mphi$ and if \eqref{main:eq1} holds for almost every $\xi\in \R$, by Lebesgue dominated convergence theorem,  then \eqref{main:eq1} holds in the sense of distributions.
If all elements in $\fPhi, \fPsi, \tilde \fPhi, \tilde \fPsi$ are essentially nonnegative measurable functions, then it is not difficult to verify that the conditions in items (i) and (ii) of Theorem~\ref{thm:main:1} are equivalent to the following simple conditions:
\begin{equation}\label{nonnegative:main:eqs:1}
\sum_{\ell=1}^{\mphi} \fphi^\ell(\df\xi)\tilde \fphi^\ell(\df\xi)+\sum_{\ell=1}^\mpsi \fpsi^\ell(\df\xi)\tilde \fpsi^\ell(\df\xi)=\sum_{\ell=1}^{\mphi} \fphi^\ell(\xi)\tilde \fphi^\ell(\xi), \qquad a.e.\; \xi\in \R
\end{equation}
and
\begin{equation}\label{nonnegative:main:eqs:2}
\begin{split}
\fphi^\ell(\xi) \tilde \fphi^\ell(\xi+2\pi k)=0 \quad &\hbox{and}\quad
\fpsi^{\ell'}(\xi)\tilde \fpsi^{\ell'}(\xi+2\pi k)=0 \qquad a.e.\; \xi\in \R, \;\\
&\forall\; k\in \Z \bs \{0\}, \; \ell=1, \ldots, \mphi,\; \ell'=1, \ldots, \mpsi.
\end{split}
\end{equation}
As we shall see in section~2, items (i) and (ii) of Theorem~\ref{thm:main:1} correspond to a natural multiresolution-like structure, which is closely linked to a fast wavelet transform. The condition in item (iii) of Theorem~\ref{thm:main:1} is a natural normalization condition which is related to \eqref{def:dwf}.

Comparing with the characterization of a pair of homogeneous dual wavelet frames in the space $\Lp{2}$ or a homogeneous orthonormal wavelet basis in $\Lp{2}$ (e.g., see
\cite{DGM,Han:msc,Han:frame,HW,RonShen:dwf,RonShen:twf}), Theorem~\ref{thm:main:1} has several interesting features. Firstly, the characterization in items (i)--(iii) of Theorem~\ref{thm:main:1} does not involve any infinite series or infinite sums; this is in sharp contrast to the homogeneous setting in $\Lp{2}$. Secondly, as we shall see in section~2, all the involved infinite sums in the proof of Theorem~\ref{thm:main:1} are in fact finite sums. This allows us to easily generalize Theorem~\ref{thm:main:1} to any real dilation factors and to nonstationary wavelets, see sections~2 and 4 for detail. Thirdly, we do not require any stability (Bessel) property of the wavelet systems, while the homogeneous setting in $\Lp{2}$ needs the stability property to guarantee the convergence of the involved infinite series. Fourthly, we do not require in Theorem~\ref{thm:main:1} that the generating wavelet functions possess any order of vanishing moments or smoothness, while all the generating wavelet functions in the homogeneous setting require at least one vanishing moment. Lastly, from a pair of nonhomogeneous dual wavelet frames in $\Lp{2}$, we shall see in section~3 that one can always derive an associated pair of homogeneous dual wavelet frames in $\Lp{2}$. In fact, most homogeneous wavelet systems in the literature are derived in such a way.
We mention that weak convergence of wavelet expansions has been characterized in \cite{FGWW} for homogeneous wavelet systems.
Similar weak convergence of wavelet expansions that are related to \eqref{def:dwf} also appeared in the study of homogeneous dual wavelet frames in $\Lp{2}$ and their frame approximation properties, for example, see \cite{DHRS,Han:msc,Han:frame}. We also point out that the approach in this paper can be extended to frequency-based homogeneous wavelet systems in the distribution space $\mathscr{D}'(\R \bs \{0\})$.

The following result generalizes the Oblique Extension Principle (OEP) and naturally connects a wavelet filter bank with a pair of frequency-based nonhomogeneous dual wavelet frames in $\Dpr$.

\begin{theorem}\label{thm:main:2} Let $\df$ be an integer such that $|\df|>1$. Let $\fa, \fth^1, \ldots, \fth^\mphi, \fb^1, \ldots, \fb^{\mpsi}$ and $\tilde \fa, \tilde \fth^1, \ldots, \tilde \fth^\mphi$, $\tilde \fb^1, \ldots, \tilde \fb^{\mpsi}$ be $2\pi$-periodic measurable functions on $\R$. Suppose that there are measurable functions $\fphi, \tilde \fphi$ satisfying
\begin{equation}\label{ref:phi}
\fphi(\df\xi)=\fa(\xi)\fphi(\xi)\quad \hbox{and}\quad \tilde \fphi(\df\xi)=\tilde \fa(\xi)\tilde \fphi(\xi), \qquad a.e.\; \xi\in \R.
\end{equation}
Define $\fPhi, \fPsi, \tilde \fPhi, \tilde \fPsi$ as in \eqref{gen:set} with
%
\begin{equation}\label{OEP:phi}
\fphi^\ell(\xi):=\fth^\ell(\xi)\fphi(\xi), \quad
\tilde \fphi^\ell(\xi):=\tilde \fth^\ell(\xi)\tilde \fphi(\xi),\qquad \ell=1, \ldots, \mphi
\end{equation}
and
\begin{equation}\label{OEP:psi}
\fpsi^\ell(\df\xi):=\fb^\ell(\xi)\fphi(\xi), \quad
\tilde \fpsi^\ell(\df \xi):=\tilde \fb^\ell(\xi)\tilde \fphi(\xi),\qquad \ell=1, \ldots, \mpsi.
\end{equation}
Assume that all the elements in $\fPhi, \fPsi, \tilde \fPhi, \tilde \fPsi$ belong to $\lLp{2}$.
Then $(\FWS_J(\fPhi; \fPsi), \FWS_J(\tilde\fPhi; \tilde \fPsi))$ forms a pair of frequency-based nonhomogeneous dual wavelet frames in the distribution space $\Dpr$ for some integer $J$ (or for all integers $J$), if and only if,
\begin{equation}\label{Phitoone}
\lim_{j\to +\infty}  \Theta(\df^{-j}\cdot)\ol{\fphi(\df^{-j}\cdot)}\tilde \fphi (\df^{-j}\cdot)=1 \quad \hbox{in the sense of distributions}
\end{equation}
with
\begin{equation}\label{Theta}
\Theta(\xi):=\sum_{\ell=1}^\mphi \ol{\fth^\ell(\xi)} \tilde \fth^\ell(\xi),
\end{equation}
and the following fundamental identities are satisfied:
\begin{equation}\label{OEP:cond:0}
\Theta(\df\xi) \ol{\fa(\xi)}  \tilde \fa(\xi)+\sum_{\ell=1}^{\mpsi} \ol{\fb^\ell(\xi)} \tilde \fb^\ell(\xi)=\Theta(\xi),
\quad a.e.\; \xi\in \sigma_\fphi\cap \sigma_{\tilde \fphi},
\end{equation}
and
\begin{equation}\label{OEP:cond:1}
\Theta(\df\xi) \ol{\fa(\xi)}  \tilde \fa(\xi+\tfrac{2\pi \omega}{\df})+\sum_{\ell=1}^{\mpsi} \ol{\fb^\ell(\xi)} \tilde \fb^\ell(\xi+\tfrac{2\pi \omega}{\df})=0,
\quad a.e.\; \xi\in \sigma_\fphi\cap (\sigma_{\tilde \fphi}-\tfrac{2\pi \omega}{\df}),
\end{equation}
for all $\omega=1, \ldots, |\df|-1$,
where $\sigma_{\tilde \fphi}-\frac{2\pi\omega}{\df}:=\{ \xi-\frac{2\pi\omega}{\df}\, : \, \xi\in \sigma_{\tilde \fphi}\}$ and
\begin{equation}\label{sigmaphi}
\sigma_\fphi:=\Big \{ \xi\in \R\; : \; \sum_{k\in \Z} |\fphi(\xi+2\pi k)|\ne 0\Big \},\quad
\sigma_{\tilde \fphi}:=\Big \{ \xi\in \R \; : \; \sum_{k\in \Z} |\tilde \fphi(\xi+2\pi k)|\ne 0\Big \}.
\end{equation}
In particular, if all $\fth^1, \ldots, \fth^\mphi, \fb^1, \ldots, \fb^\mpsi, \tilde \fth^1, \ldots, \tilde \fth^\mphi, \tilde \fb^1, \ldots, \tilde \fb^\mpsi$ are $2\pi$-periodic measurable functions in $\lLp{2}$ and if there exist positive real numbers $\tau$ and $C$ such that
\begin{equation}\label{Lip:ata}
|1-\fa(\xi)|\le C |\xi|^{\tau} \quad \hbox{and}\quad |1-\tilde \fa(\xi)|\le C |\xi|^{\tau}, \qquad a.e.\; \xi\in [-\pi, \pi],
\end{equation}
then the frequency-based standard refinable measurable functions $\fphi, \tilde \fphi$ with masks $\fa, \tilde \fa$ and the dilation factor $\df$, which are defined by
\begin{equation}\label{reffunc}
\fphi(\xi):=\prod_{j=1}^\infty \fa(\df^{-j}\xi) \quad \hbox{and}\quad \tilde \fphi(\xi):=\prod_{j=1}^\infty \tilde \fa(\df^{-j}\xi),
\end{equation}
are well-defined for almost every $\xi\in \R$ and in fact $\fphi,\tilde \fphi\in \lLp{\infty}$.
Then all elements in $\fPhi, \fPsi, \tilde \fPhi, \tilde \fPsi$ belong to $\lLp{2}$. Moreover,
$(\FWS_J(\fPhi; \fPsi), \FWS_J(\tilde \fPhi; \tilde \fPsi))$ forms a pair of frequency-based nonhomogeneous dual wavelet frames in the distribution space $\Dpr$ for some integer $J$ (or for all integers $J$), if and only if, the identities \eqref{OEP:cond:0} and \eqref{OEP:cond:1} are satisfied for all $\omega=1, \ldots, |\df|-1$, and
$\lim_{j\to +\infty}\Theta(\df^{-j}\cdot)=1$ in the sense of distributions.
\end{theorem}

Note that \eqref{Lip:ata} is automatically satisfied with $\tau=1$ if $\fa$ and $\tilde \fa$ are $2\pi$-periodic trigonometric polynomials with $\fa(0)=\tilde \fa(0)=1$.
A similar result to Theorem~\ref{thm:main:2} also holds when $\fphi=(\fphi^1, \ldots, \fphi^\mphi)^T$ and $\tilde \fphi=(\tilde \fphi^1, \ldots, \tilde \fphi^\mphi)^T$ are refinable measurable function vectors.
The identities in \eqref{OEP:cond:0} and \eqref{OEP:cond:1} with $\df=2$ and $\mphi=1$ are called the oblique extension principle in \cite{DHRS}, provided that all elements in $\fPhi, \fPsi, \tilde \fPhi, \tilde \fPsi$ belong to $\Lp{2}$ and satisfy some technical conditions to guarantee the Bessel (stability) property of the homogeneous wavelet systems $\FWS(\fPsi)$ and $\FWS(\tilde \fPsi)$ in the space $\Lp{2}$ (see \cite{CHS,DH:dwf,DHRS,Han:oep,RonShen:dwf,RonShen:twf}).
In contrast, our results here generally do not require any a priori condition on the generating wavelet functions and provide a natural explanation for the connection between the perfect reconstruction property induced by OEP in \eqref{OEP:cond:0} and \eqref{OEP:cond:1} in the discrete filter bank setting to wavelets and framelets in the function setting.

The structure of the paper is as follows. In order to prove Theorems~\ref{thm:main:1} and \ref{thm:main:2}, we shall introduce some auxiliary results in section~2. In particular, we shall provide sufficient conditions in section~2 for the absolute convergence of the infinite series in \eqref{series}. Then we shall prove
Theorems~\ref{thm:main:1} and \ref{thm:main:2} in section~2.
To explain in more detail about our motivation and importance for studying frequency-based nonhomogeneous wavelet systems, we shall discuss in section 3 nonhomogeneous wavelet systems in various function spaces such as $\Lp{2}$ and Sobolev spaces, as initiated in \cite{HanShen:ca:2009}. We shall see in section~3 that under the stability property, a pair of frequency-based nonhomogeneous dual wavelet frames can be naturally extended from the distribution space to a pair of dual function spaces. In section~3, we shall also explore the connections between nonhomogeneous and homogeneous wavelet systems in the space $\Lp{2}$. To illustrate the flexibility and generality of the approach in this paper, we further study
nonstationary wavelets which are useful in many applications, since nonstationary wavelet filter banks can be implemented in almost the same way and efficiency as a traditional fast wavelet transform. However, except a few special cases as discussed in \cite{CHS:ntwf,CohenDyn,HanShen:siam:2008}, only few theoretical results on nonstationary wavelets are available in the literature, probably partially due to the difficulty in guaranteeing the membership of the associated refinable functions in $\Lp{2}$ and in establishing the stability property of the nonstationary wavelet systems in $\Lp{2}$. In section~4, we present a complete characterization of a pair of frequency-based nonstationary dual wavelet frames in the distribution space. Though the statements and notation in section~4 on nonstationary wavelets seem a little bit more complicated comparing with the stationary case in sections 1--3,
it is worth our effort to provide a better picture to understand nonstationary wavelets, since there are few theoretical results on this topic in the literature.

To understand and study wavelet systems in various function spaces, it is our opinion that there are two key fundamental ingredients to be considered.
One ingredient is the notion investigated in this paper of a pair of frequency-based nonhomogeneous dual wavelet frames in the distribution space which enables us to completely separate its perfect reconstruction property from its stability property in function spaces.  The other ingredient is the stability issue of nonhomogeneous wavelet systems in function spaces which we didn't discuss in this paper but shall be addressed elsewhere.

\section{Frequency-based Nonhomogeneous Wavelet Systems in the Distribution Space}

In this section, we study pairs of frequency-based nonhomogeneous dual wavelet frames in the distribution space.
To prove Theorems~\ref{thm:main:1} and \ref{thm:main:2}, we first present some sufficient conditions for the absolute convergence of the infinite series in \eqref{series}.

For $1\le p\le \infty$, by $\TLp{p}$ we denote the set of all $2\pi$-periodic measurable functions $f$ such that $\int_{-\pi}^{\pi} |f(x)|^p dx<\infty$ (with the usual modification for $p=\infty$).

%
%

By the following result, we always have the absolute convergence of the infinite series in \eqref{series} provided that all the frequency-based wavelet functions are from the space $\lLp{2}$.

\begin{lemma}\label{lem:converg}
Let $\gl$ be a nonzero real number and let $\fpsi, \tilde \fpsi\in \lLp{2}$. Then for all $\ff, \fg\in \D$,
\begin{equation}\label{parseval}
\sum_{k\in \Z} \la \ff, \fpsi_{\gl;0,k} \ra \la \tilde \fpsi_{\gl;0,k}, \fg\ra =2\pi \int_\R \sum_{k\in \Z} \ff(\xi)\ol{\fg(\xi+2\pi \gl^{-1} k)}\,\ol{\fpsi(\gl\xi)} \tilde \fpsi(\gl\xi+2\pi k) d\xi
\end{equation}
with the series on the left-hand side converging absolutely. Note that the infinite sum on the right-hand side of \eqref{parseval} is in fact finite.
\end{lemma}

\begin{proof} By $L_{\infty,c}(\R)$ we denote the linear space of all compactly supported measurable functions in $\Lp{\infty}$.  Note that $\D \subseteq L_{\infty,c}(\R)$. More generally, we prove \eqref{parseval} for $\ff, \fg\in L_{\infty,c}(\R)$.
Denote
\[
\fh(\xi):=\sum_{k\in \Z} \ff(\gl^{-1}(\xi+2\pi k)) \ol{\fpsi(\xi+2\pi k)} \quad \hbox{and}\quad \tilde \fh(\xi):=\sum_{k\in \Z} \fg(\gl^{-1}(\xi+2\pi k)) \ol{\tilde \fpsi(\xi+2\pi k)}.
\]
Now we show that $\fh, \tilde \fh$ are well-defined functions in $\TLp{2}$. In fact, since $\ff\in L_{\infty,c}(\R)$,  $\ff$ has compact support and therefore, $\ff$ is essentially supported inside $[-\pi c_\ff, \pi c_\ff]$ for some $c_\ff>0$ with $c_\ff$ depending on $\ff$. Now it is easy to see that
\begin{equation}\label{H:finite}
\fh(\xi)=\sum_{k\in \Z, |k| \le (1+|\gl| c_\ff)/2}
\ff(\gl^{-1}(\xi+2\pi k)) \ol{\fpsi(\xi+2\pi k)}, \qquad \xi\in [-\pi,\pi].
\end{equation}
Since $\fpsi\in \lLp{2}$ and $\ff\in L_{\infty,c}(\R)$, we see that $\ff(\gl^{-1}(\cdot+2\pi k)) \ol{\fpsi(\cdot+2\pi k)}\in \Lp{2}$ for every $k\in \Z$. Therefore, $\fh$ is a well-defined $2\pi$-periodic function in $\TLp{2}$. Similarly, we have $\tilde \fh\in \TLp{2}$.
Note that
\[
\la \ff, \fpsi_{\gl; 0,k} \ra=|\gl|^{1/2} \int_\R \ff(\xi) \ol{\fpsi(\gl \xi)} e^{i k \gl \xi}d\xi=|\gl|^{-1/2} \int_\R \ff(\gl^{-1} \xi)\ol{\fpsi(\xi)} e^{ik\xi} d\xi=|\gl|^{-1/2} \int_{-\pi}^\pi \fh(\xi) e^{ik\xi} d\xi
\]
and $\la \fg, \tilde \fpsi_{\gl; 0, k} \ra=|\gl|^{-1/2}\int_{-\pi}^{\pi} \tilde \fh(\xi) e^{i k\xi} d\xi$. Since $\fh, \tilde \fh\in \TLp{2}$, by the Parseval identity, we have
\[
\sum_{k\in \Z} |\la \ff, \fpsi_{\gl; 0, k}\ra|^2=2\pi |\gl|^{-1} \int_{-\pi}^\pi |\fh(\xi)|^2 d\xi<\infty,\quad
\sum_{k\in \Z} |\la \fg, \tilde \fpsi_{\gl; 0, k}\ra|^2=2\pi |\gl|^{-1} \int_{-\pi}^\pi |\tilde \fh(\xi)|^2 d\xi<\infty,
\]
and
\[
\sum_{k\in \Z} \la \ff, \fpsi_{\gl; 0,k}\ra \la \tilde \fpsi_{\gl; 0, k}, \fg\ra=2\pi |\gl|^{-1} \int_{-\pi}^\pi \fh(\xi)\ol{\tilde \fh(\xi)} d\xi
\]
with the series on the left-hand side converging absolutely.
By the finite sum in \eqref{H:finite}, we have
\begin{align*}
|\gl|^{-1} \int_{-\pi}^\pi \fh(\xi)\ol{\tilde \fh(\xi)} d\xi&= |\gl|^{-1} \int_\R \ff(\gl^{-1}\xi) \ol{\fpsi(\xi)}\, \ol{\tilde \fh(\xi)} d\xi=\int_\R \ff(\xi) \ol{\fpsi(\gl \xi)}\, \ol{\tilde \fh(\gl\xi)} d\xi\\
&= \int_\R \sum_{k\in \Z} \ff(\xi)\ol{\fg(\xi+2\pi \gl^{-1} k)}\,\ol{\fpsi(\gl\xi)} \tilde \fpsi(\gl\xi+2\pi k) d\xi,
\end{align*}
which completes the proof.
\end{proof}

The condition $\fpsi, \tilde \fpsi\in \lLp{2}$ in Lemma~\ref{lem:converg} is only used in this paper to guarantee the identity \eqref{parseval}. As long as \eqref{parseval} holds for $\ff, \fg\in \D$ and the frequency-based wavelet functions belong to $\lLp{1}$, all the claims in this paper still hold.
Note that $\lLp{1}$ is just the set of all measurable functions that can be globally identified as distributions. So, $\lLp{1}$ is the most natural and weakest space for us to study wavelets and framelets in the distribution space. The condition $\fpsi, \tilde \fpsi \in \lLp{2}$ in Lemma~\ref{lem:converg} could be replaced by other conditions.
For example, for every positive integer $k$, if there exist positive numbers $\tau_k$ and $C_k$ such that $\tau_k>1/2$ and
$|\fpsi(\xi)-\fpsi(\zeta)|\le C_k |\xi-\zeta|^{\tau_k}$ for all $\xi, \zeta\in [-k,k]$,
then it is not difficult to check by \eqref{H:finite} that $\fh$ with $\ff\in \D$ is a $2\pi$-periodic Lipschitz function with some Lipschitz exponent $\tau>1/2$. By Bernstein Theorem, $\fh$ has an absolutely convergent Fourier series. Now for any $\tilde \fpsi \in \lLp{1}$, it is easy to prove that \eqref{parseval} indeed holds for all $\ff, \fg\in \D$.
Other assumptions could be used to guarantee \eqref{parseval}. But $\lLp{2}$ is a large space containing the Fourier transforms of all compactly supported distributions and of all elements in all Sobolev spaces.
For simplicity of presentation, we shall stick to the space $\lLp{2}$ for our discussion of frequency-based wavelets and framelets.

\begin{lemma}\label{lem:to1}
Let $\{\gl_j\}_{j=J}^\infty$ be a sequence of nonzero real numbers such that $\lim_{j\to +\infty} \gl_j=0$.
Let $\fphi^{j,1}, \ldots, \fphi^{j,\mphi_j}$ and $\tilde \fphi^{j,1}, \ldots, \tilde \fphi^{j,\mphi_j}$ be elements in $\lLp{2}$ with $\mphi_j\in \N$ and $j\ge J$. Then
\begin{equation}\label{lem:to1:eq1}
\lim_{j\to +\infty} \sum_{\ell=1}^{\mphi_j} \sum_{k\in \Z} \la \ff, \fphi^{j,\ell}_{\gl_j;0,k}\ra \la \tilde \fphi^{j,\ell}_{\gl_j; 0, k}, \fg\ra=2\pi \la \ff, \fg\ra \qquad \forall\, \ff, \fg\in \D,
\end{equation}
%
%
%
if and only if,
\begin{equation} \label{main:eq1:general}
\lim_{j\to +\infty} \sum_{\ell=1}^{\mphi_j} \ol{\fphi^{j,\ell}(\gl_j\cdot)} \tilde \fphi^{j,\ell}(\gl_j\cdot)=1 \qquad \hbox{in the sense of distributions}.
\end{equation}
\end{lemma}

\begin{proof} By Lemma~\ref{lem:converg}, we have
\begin{equation}\label{s2:eq6}
\sum_{\ell=1}^{\mphi_j} \sum_{k\in \Z} \la \ff, \fphi^{j,\ell}_{\gl_j;0,k} \ra \la \tilde \fphi^{j,\ell}_{\gl_j; 0,k}, \fg\ra=2\pi \int_\R \sum_{\ell=1}^{\mphi_j} \sum_{k\in \Z} \ff(\xi) \ol{ \fg(\xi+2\pi \gl_j^{-1} k)}\, \ol{\fphi^{j,\ell}(\gl_{j}\xi)}\tilde \fphi^{j,\ell}(\gl_{j}\xi+2\pi k) d\xi.
\end{equation}
By $\ff, \fg\in \D$, $\ff$ and $\fg$ are compactly supported.
Since $\lim_{j\to +\infty} \gl_j=0$, there exists an integer $J_{\ff,\fg}$ such that $\ff(\xi)\ol{\fg(\xi+2\pi \gl_{j}^{-1} k)}=0$ for all $\xi\in \R$, $k\in \Z \bs \{0\}$, and $j\ge J_{\ff,\fg}$.
That is, for $j\ge \max(J,J_{\ff,\fg})$, \eqref{s2:eq6} becomes
\begin{equation}\label{s2:eq7}
\sum_{\ell=1}^{\mphi_j} \sum_{k\in \Z} \la \ff, \fphi^{j,\ell}_{\gl_j; 0,k}\ra \la \tilde \fphi^{j,\ell}_{\gl_j;0,k}, \fg\ra
=2\pi \int_\R \ff(\xi)\ol{\fg(\xi)} \sum_{\ell=1}^{\mphi_j} \ol{\fphi^{j,\ell}(\gl_{j}\xi)} \tilde \fphi^{j,\ell}(\gl_j \xi) d\xi.
\end{equation}
If \eqref{main:eq1:general} holds in the sense of distributions, then it follows directly from \eqref{s2:eq7} that \eqref{lem:to1:eq1} holds.

Conversely, we can take $\fg\in \D$ such that $\fg$ takes value one on the support of $\ff$, now it follows from \eqref{lem:to1:eq1} and \eqref{s2:eq7} that
\[
\lim_{j\to +\infty} \int_\R \ff(\xi) \sum_{\ell=1}^{\mphi_j} \ol{\fphi^{j,\ell}(\gl_j\xi)} \tilde \fphi^{j,\ell}(\gl_j \xi) d\xi
=\lim_{j\to+\infty}\frac{1}{2\pi} \sum_{\ell=1}^{\mphi_j} \sum_{k\in \Z} \la \ff, \fphi^{j,\ell}_{\gl_j;0,k}\ra \la \tilde \fphi^{j,\ell}_{\gl_j;0,k}, \fg\ra
=\la \ff, \fg\ra=\la \ff,1\ra.
\]
Hence, \eqref{main:eq1:general} holds in the sense of distributions.
\end{proof}

In the next auxiliary result, we shall study a multiresolution-like structure. More precisely, we have the following result.

\begin{lemma}\label{lem:twolevel}
Let $\gl$ be a nonzero real number. Let $\fphi^1, \ldots, \fphi^{\mphi}$, $\fpsi^1, \ldots, \fpsi^{\mpsi}$, $\feta^1, \ldots, \feta^{\meta}$ and $\tilde \fphi^1, \ldots, \tilde \fphi^{\mphi}$, $\tilde \fpsi^1, \ldots, \tilde \fpsi^{\mpsi}, \tilde \feta^1, \ldots, \tilde \feta^{\meta}$ be elements in $\lLp{2}$.
Then
\begin{equation}\label{twolevel:general}
\begin{split}
\sum_{\ell=1}^\mphi \sum_{k\in \Z} \la \ff, \fphi^{\ell}_{1;0,k}\ra\la \tilde \fphi^{\ell}_{1;0,k}, \fg\ra
&+\sum_{\ell=1}^\mpsi \sum_{k\in \Z} \la \ff, \fpsi^{\ell}_{1;0,k}\ra\la \tilde \fpsi^{\ell}_{1;0,k}, \fg\ra\\
&=\sum_{\ell=1}^\meta \sum_{k\in \Z} \la \ff, \feta^{\ell}_{\gl;0,k}\ra\la \tilde \feta^{\ell}_{\gl;0,k}, \fg\ra \qquad \forall\; \ff, \fg\in \D
\end{split}
\end{equation}
if and only if
\begin{align}
&I_\fphi^{k}(\xi)+I^{k}_{\fpsi}( \xi)=I^{\gl k}_{\feta}(\gl \xi), \qquad a.e.\, \xi\in \R,\; \forall\; k\in \Z \cap [\gl^{-1} \Z], \label{I:eq1}\\
&I_\fphi^{k}(\xi)+I^{k}_{\fpsi}(\xi)=0, \qquad a.e.\, \xi\in \R,\; \forall\; k\in \Z \bs [\gl^{-1} \Z], \label{I:eq2}\\
&I^{\gl k}_\feta(\xi)=0, \qquad a.e.\, \xi\in \R,\; \forall \; k\in [\gl^{-1} \Z] \bs \Z, \label{I:eq3}
\end{align}
where $\gl^{-1} \Z:=\{ \gl^{-1} k\; : \; k\in \Z\}$ and
$I^{\gl k}_\feta(\xi):=\sum_{\ell=1}^\meta \ol{\feta^\ell (\xi)} \tilde \feta^\ell(\xi+2\pi \gl k)$, $k\in \gl^{-1}\Z$, and
\begin{equation}\label{I}
I^k_\fphi(\xi):=\sum_{\ell=1}^\mphi \ol{\fphi^\ell (\xi)} \tilde \fphi^\ell(\xi+2\pi k), \qquad
I^k_\fpsi(\xi):=\sum_{\ell=1}^\mpsi \ol{\fpsi^\ell (\xi)} \tilde \fpsi^\ell(\xi+2\pi k), \qquad k\in \Z.
\end{equation}
\end{lemma}

\begin{proof} By Lemma~\ref{lem:converg}, all the infinite series in \eqref{twolevel:general} converge absolutely and \eqref{twolevel:general} is equivalent to
\[
\int_\R \sum_{k\in \Z} \ff(\xi) \ol{\fg(\xi+2\pi k)} \Big(I^{k}_\fphi(\xi)
+I^{k}_\fpsi(\xi)\Big) d\xi=\int_\R \sum_{k\in [\gl^{-1} \Z]} \ff(\xi) \ol{\fg(\xi+2\pi k)} I^{\gl k}_\feta(\gl \xi) d\xi,
\]
which can be easily rewritten as
\begin{equation}\label{s2:eq2}
\begin{split}
&\int_\R \sum_{k\in \Z \cap [\gl^{-1} \Z]} \ff(\xi) \ol{\fg(\xi+2\pi k)} \Big(I^{k}_\fphi(\xi)
+I^{k}_\fpsi(\xi)-I^{\gl k}_\feta(\gl \xi)
\Big) d\xi\\
&+\int_\R \sum_{k\in \Z \bs [\gl^{-1} \Z]} \ff(\xi) \ol{\fg(\xi+2\pi k)} \Big(I^{k}_\fphi(\xi)
+I^{k}_\fpsi(\xi)\Big) d\xi
=\int_\R \sum_{k\in [\gl^{-1} \Z]\bs \Z } \ff(\xi) \ol{\fg(\xi+2\pi k)} I^{\gl k}_\feta(\gl \xi) d\xi.
\end{split}
\end{equation}

Sufficiency.
If \eqref{I:eq1}, \eqref{I:eq2}, and \eqref{I:eq3} are  satisfied, then it is obvious that \eqref{s2:eq2} is true and therefore, \eqref{twolevel:general} holds.

Necessity. Denote $\gL:= \Z \cup [\gl^{-1} \Z]$. For a point $x\in \R$, we define $\hbox{dist}(x, \gL):=\inf_{y\in \gL} |x-y|$.
By \eqref{twolevel:general}, \eqref{s2:eq2} holds.
Let $k_0\in \Z \cap [\gl^{-1} \Z]$ and $\xi_0\in \R$ be temporarily fixed. Then it is easy to check that
$\gep:=\frac{\pi}{2}\hbox{dist}(k_0, \gL \bs \{k_0\})>0$. Consider all $\ff,\fg\in \D$ such that the support of $\ff$ is contained inside $(\xi_0-\gep, \xi_0+\gep)$ and the support of $\fg$ is contained inside $(\xi_0-2\pi k_0-\gep, \xi_0-2\pi k_0+\gep)$.
Then it is not difficult to verify that
\[
\ff(\xi) \ol{\fg(\xi+2\pi k)}=0\quad  \forall\; \xi\in \R, k\in \gL \bs \{k_0\},
\]
from which we see that \eqref{s2:eq2} becomes
\begin{equation}\label{s2:eq3}
\int_\R \ff(\xi)\ol{\fg(\xi+2\pi k_0)} (I^{k_0}_\fphi(\xi)
+I^{k_0}_\fpsi(\xi)-I^{\gl k_0}_\feta(\gl \xi)) d\xi=0
\end{equation}
for all $\ff, \fg\in \D$ such that $\mbox{supp}\,\ff\subseteq (\xi_0-\gep, \xi_0+\gep)$ and
$\mbox{supp}\,\fg\subseteq (\xi_0-2\pi k_0-\gep, \xi_0-2\pi k_0+\gep)$. From \eqref{s2:eq3}, we must have $I^{k_0}_\fphi(\xi)+I^{k_0}_\fpsi(\xi)
-I^{\gl k_0}_\feta(\gl \xi)=0$ for almost every $\xi\in (\xi_0-\gep, \xi_0+\gep)$. Thus, \eqref{I:eq1} must be true. \eqref{I:eq2} and \eqref{I:eq3} can be proved by the same argument.
\end{proof}

Now we have the generalized version of Theorem~\ref{thm:main:1} with a general real dilation factor $\df$.

\begin{theorem}\label{thm:main:1:general}
Let $\df$ be a real number such that $|\df|>1$.
Let $\fPhi, \fPsi, \tilde \fPhi, \tilde \fPsi$ in \eqref{gen:set} be subsets of $\lLp{2}$. Then
$(\FWS_{J}(\Phi; \Psi), \FWS_{J}(\tilde \Phi; \tilde \Psi))$, where
$\FWS_J(\Phi; \Psi)$ is defined in \eqref{FWS},
forms a pair of frequency-based nonhomogeneous dual wavelet frames in the distribution space $\Dpr$ for some integer $J$ (or for all integers $J$), if and only if,
\begin{equation}\label{main:eq1:dm}
\lim_{j\to +\infty} \sum_{\ell=1}^r \ol{\fphi^\ell(\df^{-j}\cdot)} \tilde \fphi^\ell(\df^{-j}\cdot)=1 \qquad \hbox{in the sense of distributions}
\end{equation}
and
\begin{equation}\label{I:eq1:sp}
\begin{split}
\sum_{\ell=1}^\mphi \ol{\fphi^\ell(\xi)} \tilde \fphi^\ell(\xi+2\pi k)
&+\sum_{\ell=1}^\mpsi \ol{\fpsi^\ell(\xi)} \tilde \fpsi^\ell(\xi+2\pi k)\\
&=\sum_{\ell=1}^\mphi \ol{\fphi^\ell(\df^{-1}\xi)} \tilde \fphi^\ell(\df^{-1}(\xi+2\pi k)),
\quad a.e.\, \xi\in \R,\; \forall\; k\in \Z \cap [\df \Z],
\end{split}
\end{equation}
\begin{align}
&\sum_{\ell=1}^\mphi \ol{\fphi^\ell(\xi)} \tilde \fphi^\ell(\xi+2\pi k)
+\sum_{\ell=1}^\mpsi \ol{\fpsi^\ell(\xi)} \tilde \fpsi^\ell(\xi+2\pi k)=0, \qquad a.e.\, \xi\in \R,\; \forall\; k\in \Z \bs [\df \Z], \label{I:eq2:sp}\\
&\sum_{\ell=1}^\mphi \ol{\fphi^\ell(\xi)} \tilde \fphi^\ell(\xi+2\pi \df^{-1} k)=0, \qquad a.e.\, \xi\in \R,\; \forall \; k\in [\df \Z] \bs \Z. \label{I:eq3:sp}
\end{align}
\end{theorem}

\begin{proof} By the following simple observation, we have
\begin{equation}\label{rel:dilation}
\la \ff_{\gl; n, k}, \fg_{\gl; n,k}\ra=\la \ff, \fg\ra \quad \mbox{and}\quad \la \ff_{\gl; 0, 0}, \fpsi^\ell_{\gl';n,k}\ra
=\la \ff, \fpsi^\ell_{\gl'\gl^{-1}; n,k}\ra,\qquad \gl, \gl' \in \R \bs \{0\}, k,n \in \R.
\end{equation}
Now it is straightforward to see that for all $\ff, \fg\in \D$,
\begin{equation}\label{MRA}
\begin{split}
\sum_{\ell=1}^\mphi \sum_{k\in \Z} \la \ff, \fphi^{\ell}_{\df^{-j};0,k}\ra \la \tilde \fphi^{\ell}_{\df^{-j};0,k}, \fg\ra
&+\sum_{\ell=1}^\mpsi \sum_{k\in \Z} \la \ff, \fpsi^{\ell}_{\df^{-j};0,k}\ra\la \tilde \fpsi^{\ell}_{\df^{-j};0,k}, \fg\ra\\
&= \sum_{\ell=1}^\mphi \sum_{k\in \Z} \la \ff, \fphi^{\ell}_{\df^{-j-1};0,k}\ra\la \tilde \fphi^{\ell}_{\df^{-j-1};0,k}, \fg\ra,
\end{split}
\end{equation}
if and only if, for all $\ff, \fg\in \D$,
\begin{equation}\label{twolevel}
\sum_{\ell=1}^\mphi \sum_{k\in \Z} \la \ff, \fphi^{\ell}_{1;0,k}\ra \la \tilde \fphi^{\ell}_{1;0,k}, \fg\ra
+\sum_{\ell=1}^\mpsi \sum_{k\in \Z} \la \ff, \fpsi^{\ell}_{1;0,k}\ra\la \tilde \fpsi^{\ell}_{1;0,k}, \fg\ra
= \sum_{\ell=1}^\mphi \sum_{k\in \Z} \la \ff, \fphi^{\ell}_{\df^{-1};0,k}\ra\la \tilde \fphi^{\ell}_{\df^{-1};0,k}, \fg\ra.
\end{equation}

Sufficiency. By Lemma~\ref{lem:twolevel} with $\gl=\df^{-1}$, we see that \eqref{twolevel} holds and therefore, \eqref{MRA} holds for all $\ff, \fg\in \D$ and $j\in \Z$. For $J'>J$, we define
\begin{equation}\label{SFG}
S_J^{J'}(\ff,\fg):=\sum_{\ell=1}^{\mphi} \sum_{k\in \Z} \la \ff, \fphi^{\ell}_{\df^{-J};0,k}\ra \la \tilde \fphi^{\ell}_{\df^{-J};0,k}, \fg\ra+
\sum_{j=J}^{J'-1} \sum_{\ell=1}^{\mpsi}\sum_{k\in \Z}
\la \ff, \fpsi^{\ell}_{\df^{-j};0,k}\ra \la \tilde \fpsi^{\ell}_{\df^{-j};0,k}, \fg\ra.
\end{equation}
Now by \eqref{MRA}, we can easily deduce that
\begin{equation}\label{s2:eq5}
S_J^{J'}(\ff,\fg)=\sum_{\ell=1}^{\mphi} \sum_{k\in \Z} \la \ff, \fphi^{\ell}_{\df^{-J'};0,k}\ra \la \tilde \fphi^{\ell}_{\df^{-J'};0,k}, \fg\ra.
\end{equation}
By Lemma~\ref{lem:to1}, it follows from \eqref{main:eq1:dm} that
$\lim_{J'\to +\infty} S_J^{J'}(\ff,\fg)
=2\pi \la \ff, \fg\ra$ for all $\ff, \fg\in \D$.
Hence, $(\FWS_{J}(\fPhi; \fPsi), \FWS_{J}(\tilde \fPhi; \tilde \fPsi))$ forms a pair of frequency-based nonhomogeneous dual wavelet frames in the distribution space $\Dpr$.

Necessity. By \eqref{rel:dilation},
we can easily deduce that $(\FWS_J(\fPhi; \fPsi), \FWS_J(\tilde \fPhi; \tilde \fPsi))$ forms a pair of frequency-based nonhomogeneous dual wavelet frames in the distribution space $\Dpr$ for some integer $J$ if and only if it is true for all integers $J$. Considering the difference between two consecutive integers $J$ and $J+1$, we see that \eqref{MRA} must hold and therefore, \eqref{twolevel} holds. Now by Lemma~\ref{lem:twolevel}, \eqref{I:eq1}, \eqref{I:eq2}, and \eqref{I:eq3} hold with $\feta=\fphi$ and $\meta=\mphi$,
or equivalently, \eqref{I:eq1:sp}, \eqref{I:eq2:sp}, and \eqref{I:eq3:sp} hold. Since \eqref{MRA} holds, we deduce that \eqref{s2:eq5} holds.
By Lemma~\ref{lem:to1}, it follows from our assumption $\lim_{J'\to +\infty} S_{J}^{J'}(\ff,\fg)=2\pi \la \ff, \fg\ra$ that \eqref{main:eq1:dm} must hold.
\end{proof}

We point out that Theorem~\ref{thm:main:1:general} and the approach in this paper can be extended to frequency-based homogeneous wavelet systems in the distribution space $\mathscr{D}'(\R \bs \{0\})$, which is the dual space of the test function space $\mathscr{D}(\R \bs \{0\})$ consisting of all compactly supported $C^\infty$ functions whose supports are contained inside $\R \bs \{0\}$. We shall address this issue elsewhere.

Now we are ready to prove Theorems~\ref{thm:main:1} and \ref{thm:main:2}.

\begin{proof}[Proof of Theorem~\ref{thm:main:1}] By Theorem~\ref{thm:main:1:general} with $\df$ being an integer, it suffices to show that \eqref{main:eq2} and \eqref{main:eq3} are equivalent to the three conditions \eqref{I:eq1:sp}, \eqref{I:eq2:sp}, and \eqref{I:eq3:sp}. Note that $[\df \Z] \bs \Z$ is the empty set. So, \eqref{I:eq3:sp} with $\df$ being an integer is automatically true.
It is evident that \eqref{I:eq2:sp} is equivalent to \eqref{main:eq3}. Since $\df$ is an integer, we have $\Z \cap [\df\Z]=\df \Z$. Now it is also easy to see that \eqref{main:eq2} is equivalent to \eqref{I:eq1:sp}. This completes the proof.
\end{proof}

We use Theorem~\ref{thm:main:1} to prove Theorem~\ref{thm:main:2} as follows:

\begin{proof}[Proof of Theorem~\ref{thm:main:2}] We prove the first part of Theorem~\ref{thm:main:2} first. By \eqref{ref:phi} and \eqref{OEP:phi}, we have
\[
\sum_{\ell=1}^\mphi \ol{\fphi^\ell(\df\xi)} \tilde \fphi^\ell(\df(\xi+2\pi k))
=\sum_{\ell=1}^\mphi \ol{\fth^\ell(\df\xi)} \tilde \fth^\ell(\df \xi) \ol{\fphi(\df\xi)}\tilde \fphi(\df(\xi+2\pi k))=
\ol{\fphi(\xi)}\tilde \fphi(\xi+2\pi k) \Theta(\df\xi) \ol{ \fa(\xi)}  \tilde \fa(\xi)
\]
and similarly by \eqref{OEP:psi},
\[
\sum_{\ell=1}^\mpsi \ol{\fpsi^\ell(\df\xi)}\tilde \fpsi^\ell(\df(\xi+2\pi k))=\ol{\fphi(\xi)}\tilde {\fphi}(\xi+2\pi k)\sum_{\ell=1}^\mpsi \ol{\fb^\ell(\xi)} \tilde \fb^\ell(\xi)
\]
for all integers $k\in \Z$.
Now \eqref{main:eq2} is equivalent to
\begin{equation}\label{s2:eq9}
\ol{\fphi(\xi)}\tilde {\fphi}(\xi+2\pi k)\Big(
\Theta(\df\xi) \ol{\fa(\xi)}  \tilde \fa(\xi)+ \sum_{\ell=1}^\mpsi \ol{\fb^\ell(\xi)} \tilde \fb^\ell(\xi)\Big)
=\ol{\fphi(\xi)}\tilde {\fphi}(\xi+2\pi k) \Theta(\xi)
\end{equation}
for all $k\in \Z$.
Now it is not difficult to deduce that \eqref{s2:eq9} is equivalent to \eqref{OEP:cond:0}.

Note that any $k_0\in \Z \bs [\df \Z]$ can be uniquely written as $k_0=\omega+\df k$ for $\omega\in \{1, \ldots, |\df|-1\}$ and $k\in \Z$.
Replacing $\xi$ in \eqref{main:eq3} by $\df \xi$, by the same argument, we see that \eqref{main:eq3} is equivalent to
\begin{equation}\label{s2:eq10}
\ol{\fphi(\xi)}\tilde {\fphi}(\xi+\tfrac{2\pi \omega}{\df}+2\pi k)\Big(
\Theta(\df\xi) \ol{\fa(\xi)}  \tilde \fa(\xi+\tfrac{2\pi\omega}{\df})+ \sum_{\ell=1}^\mpsi \ol{\fb^\ell(\xi)} \tilde \fb^\ell(\xi+\tfrac{2\pi \omega}{\df})\Big)
=0
\end{equation}
for all $k\in \Z$ and $\omega=1, \ldots, |\df|-1$.
Now it is not difficult to deduce that \eqref{s2:eq10} is equivalent to \eqref{OEP:cond:1} with $\omega=1, \ldots, |\df|-1$.

We now prove the second part of Theorem~\ref{thm:main:2}.
Note that \eqref{Lip:ata} implies that $\fa, \tilde \fa \in \Lp{\infty}$ and $\max(|1-\fa(\df^{-j}\xi)|, |1-\tilde \fa(\df ^{-j}\xi)|) \le C |\df|^{-\tau j}|\xi|^\tau$ for all $\xi\in (-|\df|^j\pi, |\df|^j\pi]$. It is a standard argument to show that both $\fphi$ and $\tilde \fphi$ in \eqref{ref:phi} are well-defined measurable functions in $\lLp{\infty}$ and
\begin{equation}\label{s2:eq0}
\lim_{j\to +\infty} \fphi(\df^{-j}\xi)=1=\lim_{j\to +\infty} \tilde \fphi(\df^{-j}\xi),\qquad a.e.\, \xi\in \R.
\end{equation}
In fact, by the same argument as in \cite[Page 93]{Han:siam:2008} or \cite[Page 932]{HanShen:siam:2008}, \eqref{Lip:ata} also implies that for any $\gep>0$, there exists $c>0$ such that %
\begin{equation}\label{s2:eq11}
1-\gep \le |\fphi(\xi)|\le 1+\gep \quad \hbox{and}\quad
1-\gep \le |\tilde \fphi(\xi)|\le 1+\gep, \qquad a.e.\; \xi\in [-c,c].
\end{equation}
See the proof of Theorem~\ref{thm:mra:ndwf} in section~4 for more detail on proving \eqref{s2:eq0} and \eqref{s2:eq11}.
Take $\gep=1/2$ in \eqref{s2:eq11}. By the definition of $\Theta$ in \eqref{Theta},
we conclude that
\begin{equation}\label{Theta:inq}
|\Theta(\df^{-j}\xi)\ol{\fphi(\df^{-j}\xi)}\tilde \fphi(\df^{-j}\xi)|
\le \frac{9}{4} \big |\Theta(\df^{-j}\xi)\big|
\le 9\big|\Theta(\df^{-j}\xi)\ol{\fphi(\df^{-j}\xi)}\tilde \fphi(\df^{-j}\xi)\big|
\end{equation}
for almost every $\xi\in [-c, c]$ and for all $j\ge 0$.
Consequently, by the generalized Lebesgue dominated convergence theorem and using \eqref{s2:eq0} and \eqref{Theta:inq},
we can conclude that \eqref{Phitoone} holds in the sense of distributions if and only if $\lim_{j\to +\infty} \Theta(\df^{-j}\cdot)=1$ holds in the sense of distributions.

Since $\fphi, \tilde \fphi\in \lLp{\infty}$ and $\fth^1,\ldots, \fth^\mphi, \fb^1, \ldots, \fb^{\mpsi}, \tilde \fth^1,\ldots, \tilde \fth^\mphi, \tilde \fb^1, \ldots, \tilde \fb^\mpsi \in \lLp{2}$,
by the definition in \eqref{OEP:phi} and \eqref{OEP:psi},
it is evident that all the measurable functions in $\fPhi, \fPsi, \tilde \fPhi, \tilde \fPsi$ belong to the desired space $\lLp{2}$. By what has been proved for the first part of Theorem~\ref{thm:main:2}, the claim in the second part of Theorem~\ref{thm:main:2} holds.
\end{proof}

\section{Nonhomogeneous Wavelet Systems in Function Spaces and $\Lp{2}$}

In this section, we shall discuss nonhomogeneous wavelet systems in the space $\Lp{2}$ and other function spaces. We shall see that a pair of frequency-based nonhomogeneous dual wavelet frames in the distribution space plays a basic role in the study of wavelets and framelets in various function spaces. We shall also discuss the connections between nonhomogeneous and homogeneous wavelet systems in $\Lp{2}$.

Let us recall some necessary definitions first. For $\tau\in \R$, we denote by $\HH{\tau}$ the Sobolev space consisting of all tempered distributions $f$ such that
\begin{equation}\label{sobolev}
\|f\|^2_{\HH{\tau}}:=\frac{1}{2\pi} \int_\R |\hat f(\xi)|^2 (1+|\xi|^2)^\tau d\xi<\infty.
\end{equation}
Note that $H^{-\tau}(\R)$ is the dual space of $H^\tau(\R)$, since $f\in H^{-\tau}(\R)$ can be regarded as a continuous linear functional on $H^\tau(\R)$ in the sense of $\la f, g\ra=\frac{1}{2\pi} \int_\R \hat f(\xi) \ol{\hat g(\xi)}d\xi$ for $g\in H^\tau(\R)$.

Let $\df$ be a nonzero real number and $\WS_J(\Phi; \Psi)$ be a nonhomogeneous wavelet system defined in \eqref{WS}.
We define a normed sequence space $\seqspace_{\HH{\tau}}$, indexed by the elements of $\WS_{J}(\Phi; \Psi)$, with weighted  norm as follows:
\begin{equation}\label{seq:weight:norm}
\left \| \{ w_h \}_{h\in \WS_{J}(\Phi; \Psi)} \right \|_{\seqspace_{\HH{\tau}}}^2:=
\sum_{\phi\in \Phi} \sum_{k\in \Z} |\df|^{-2\tau J} |w_{\phi_{\df^{J}; k}}|^2
+\sum_{j=J}^\infty \sum_{\psi\in \Psi} \sum_{k\in \Z} |\df|^{-2\tau j}|w_{\psi_{\df^j; k}}|^2.
\end{equation}
For $\Phi, \Psi \subseteq \HH{\tau}$, we say that $\WS_J(\Phi; \Psi)$ has {\it stability in $\HH{\tau}$} with respect to the normed sequence space $\seqspace_{\HH{\tau}}$
if there exist positive constants $C_1$ and $C_2$ such that
\begin{equation}\label{sobolevframe:0}
C_1\|f\|^2_{\HH{-\tau}}\le \| \{\la f, h\ra\}_{h\in \WS_J(\Phi;\Psi)}\|^2_{\seqspace_{\HH{\tau}}}
\le C_2 \|f\|^2_{\HH{-\tau}} \qquad \forall\; f\in \HH{-\tau},
\end{equation}
or more explicitly,
\begin{equation}\label{sobolevframe:1}
C_1\|f\|^2_{\HH{-\tau}}\le \sum_{\phi\in \Phi} \sum_{k\in \Z} |\df|^{-2\tau J}|\la f, \phi_{\df^J;k}\ra|^2+\sum_{j=J}^\infty \sum_{\psi\in \Psi}  \sum_{k\in \Z} |\df|^{-2\tau j}|\la f, \psi_{\df^j; k}\ra|^2\le C_2 \|f\|^2_{\HH{-\tau}}.
\end{equation}
Note that $\HH{\tau}$ is a Hilbert space under the inner product
\[
\la f, g\ra_{\HH{\tau}}:=\frac{1}{2\pi} \int_\R \hat f(\xi) \ol{\hat g(\xi)} (1+|\xi|^2)^\tau d\xi, \qquad f, g\in \HH{\tau}.
\]
It was shown in \cite[Proposition~2.1]{HanShen:ca:2009} that \eqref{sobolevframe:1} holds for all $f\in \HH{-\tau}$, if and only if,
\begin{equation}\label{sobolevframe:2}
\begin{split}
C_1 \|g\|_{\HH{\tau}}^2\le \sum_{\phi\in \Phi} \sum_{k\in \Z} |\la g, |\df|^{-\tau J}\phi_{\df^J;k}\ra_{\HH{\tau}}|^2&+\sum_{j=J}^\infty \sum_{\psi\in \Psi} \sum_{k\in \Z} |\la g, |\df|^{-\tau j}\psi_{\df^j;k}\ra_{\HH{\tau}}|^2\\
&\le C_2 \|g\|_{\HH{\tau}}^2,\qquad g\in \HH{\tau}.
\end{split}
\end{equation}
In other words, that $\WS_{J}(\Phi; \Psi)$ has stability in $\HH{\tau}$ with respect to $\seqspace_{\HH{\tau}}$ is equivalent to saying that after a proper renormalization of $\WS_{J}(\Phi; \Psi)$,
\begin{equation}\label{WS:sobolev}
\WS^{\tau}_J(\Phi; \Psi):=
\{|\df|^{-\tau J} \phi_{\df^J;k}\; : \; k\in \Z, \phi\in \Phi\}
\cup\{ |\df|^{-\tau j}\psi_{\df^j;k}\; : \; j\ge J, k\in \Z, \psi\in \Psi\}
\end{equation}
is a frame in the Hilbert space $\HH{\tau}$.
For $\psi\in \HH{\tau}$ and $|\df|>1$, in the following we show that
\begin{equation}\label{psi:norm}
A_1 \|\psi\|_{\HH{\tau}}^2 \le \| |\df|^{-\tau j} \psi_{\df^j;k}\|_{\HH{\tau}}^2\le A_2 \|\psi\|_{\HH{\tau}}^2 \qquad \forall\; j\ge J, k\in \Z,
\end{equation}
for some positive constants $A_1$ and $A_2$ depending on $\df, \tau, J$ and $\psi$. For $\tau<0$, to prove \eqref{psi:norm}, we further assume that
\begin{equation}\label{psi:vm}
\int_\R |\hat \psi(\xi)|^2 |\xi|^{2\tau}d\xi<\infty,
\end{equation}
which is also a necessary condition for \eqref{psi:norm} to hold. We now prove \eqref{psi:norm}. In fact,
\begin{equation}\label{psijk:norm}
\| |\df|^{-\tau j} \psi_{\df^j;k}\|_{\HH{\tau}}^2
=\frac{1}{2\pi} \int_\R |\hat \psi(\xi)|^2 (1+|\xi|^2)^\tau \Big(\frac{|\df|^{-2j}+|\xi|^2}{1+|\xi|^2}\Big)^\tau d\xi.
\end{equation}
Since $|\df|>1$, it is easy to deduce that
\begin{equation}\label{psijk:inq}
\frac{|\xi|^2}{1+|\xi|^2} \le \frac{|\df|^{-2j}+|\xi|^2}{1+|\xi|^2}
\le |\df|^{-2J}+1, \qquad \forall\; j\ge J, \xi\in \R.
\end{equation}
For $\tau\ge 0$, noting that there exists $\gep>0$ such that
\[
\frac{1}{2\pi} \int_{ \{ \zeta\in \R\, : \, |\zeta|\ge \gep \}} |\hat \psi(\xi)|^2 (1+|\xi|^2)^\tau d\xi\ge \frac{1}{2} \|\psi\|_{\HH{\tau}}^2,
\]
we deduce from \eqref{psijk:norm} and \eqref{psijk:inq} that \eqref{psi:norm} holds with $A_1=(1+\gep^{-2})^{-\tau}/2$ and $A_2=(|\df|^{-2J}+1)^\tau$.

For $\tau<0$, by our assumption in \eqref{psi:vm}, there exists $\gep>0$ such that
\[
\frac{1}{2\pi} \int_{\{ \zeta\in \R\, : \, |\zeta|\le \gep\} } |\hat \psi(\xi)|^2 |\xi|^{2\tau} d\xi\le \frac{1}{2} \|\psi\|_{\HH{\tau}}^2.
\]
Now we can easily deduce from \eqref{psijk:norm} and \eqref{psijk:inq} that \eqref{psi:norm} holds with $A_1=(|\df|^{-2J}+1)^\tau$ and $A_2=\frac{1}{2}+(1+\gep^{-2})^{-\tau}$.
Using \eqref{psijk:norm} and applying Fatou's lemma to \eqref{psi:norm} with $j\to +\infty$, it is easy to see that \eqref{psi:vm} is a necessary condition for \eqref{psi:norm} to hold.

Hence, the weights $|\df|^{-\tau j}$ in the definition of $\|\cdot\|_{\seqspace_{\HH{\tau}}}$ is chosen to be the normalization factors of the elements $\psi_{\df^j;k}$ in the space $\HH{\tau}$.
The stability property in \eqref{sobolevframe:1}, which is equivalent to the frame property in \eqref{sobolevframe:2}, characterizes the Sobolev space $\HH{-\tau}$.
We point out that the characterization of homogeneous Besov spaces by homogeneous framelets in \cite{BGN,E3} uses nonlinear approximation and is quite different in nature to \eqref{sobolevframe:1} and \eqref{sobolevframe:2}.

Now we are ready to recall the definition in \cite{HanShen:ca:2009} on a pair of nonhomogeneous dual wavelet frames in a pair of dual Sobolev spaces $(\HH{\tau}, \HH{-\tau})$. Let $\df$ be the dilation factor. Let
\begin{equation}\label{Phi:Psi}
\Phi=\{\phi^1, \ldots, \phi^\mphi\}, \quad
\Psi=\{\psi^1, \ldots, \psi^\mpsi\} \quad \mbox{and}\quad
\tilde \Phi=\{\tilde \phi^1, \ldots, \tilde \phi^\mphi\}, \quad
\tilde \Psi=\{\tilde\psi^1, \ldots, \tilde \psi^\mpsi\}
\end{equation}
be subsets of tempered distributions. We say that
the pair $(\WS_{J}(\Phi; \Psi), \WS_{J}(\tilde \Phi; \tilde \Psi))$,
or more precisely $(\WS^{\tau}_{J}(\Phi; \Psi), \WS^{ -\tau}_{J}(\tilde \Phi; \tilde \Psi))$,
forms {\it a pair of nonhomogeneous dual wavelet frames in a pair of dual Sobolev spaces $(\HH{\tau}, \HH{-\tau})$} if
\begin{enumerate}
\item $\Phi, \Psi\subseteq \HH{\tau}$ and $\WS_{J}(\Phi;\Psi)$ has stability in $\HH{\tau}$ with respect to $\seqspace_{\HH{\tau}}$. That is, $\WS^{\tau}_{J}(\Phi; \Psi)$ is a frame in $\HH{\tau}$;
\item $\tilde \Phi, \tilde \Psi \subseteq \HH{-\tau}$ and $\WS_{J}(\tilde \Phi; \tilde \Psi)$ has stability in $\HH{-\tau}$ with respect to $\seqspace_{\HH{-\tau}}$. That is, $\WS^{-\tau}_{J}(\tilde \Phi; \tilde \Psi)$ is a frame in $\HH{-\tau}$;
\item for all $f \in \HH{-\tau}$ and $g\in \HH{\tau}$, the following identity holds
\begin{equation}\label{identity}
\la f, g\ra=\sum_{\ell=1}^\mphi \sum_{k\in \Z} \la f, \phi^\ell_{\df^J;k}\ra\la \tilde \phi^\ell_{\df^J;k}, g\ra+\sum_{j=J}^\infty \sum_{\ell=1}^\mpsi \sum_{k\in \Z}
\la f, \psi^\ell_{\df^j;k}\ra \la \tilde \psi^\ell_{\df^j;k}, g\ra
\end{equation}
with the series on the right-hand side converging absolutely.
\end{enumerate}

The above definition is introduced in \cite{HanShen:ca:2009}. 
When $\tau=0$, since $\Lp{2}=\HH{0}$, it is easy to see that the above definition of a pair of nonhomogeneous dual wavelet frames in $(\Lp{2}, \Lp{2})$ becomes the definition of a pair of nonhomogeneous dual wavelet frames in $\Lp{2}$.

Suppose that the pair $(\WS_{J}(\Phi; \Psi), \WS_{J}(\tilde \Phi; \tilde \Psi))$ forms a pair of nonhomogeneous dual wavelet frames in $(\HH{\tau}, \HH{-\tau})$. Using \eqref{identity} and the upper stability (that is, the right-hand inequality in \eqref{sobolevframe:1} holds) of the two nonhomogeneous wavelet systems, it is not difficult to see that we have the following representations:
\[
f=\sum_{\ell=1}^\mphi \sum_{k\in \Z} \la f, \phi^\ell_{\df^J;k}\ra \tilde \phi^\ell_{\df^J;k}+\sum_{j=J}^\infty  \sum_{\ell=1}^\mpsi \sum_{k\in \Z}
\la f, \psi^\ell_{\df^j;k}\ra \tilde \psi^\ell_{\df^j;k}, \qquad f\in \HH{-\tau}
\]
with the series converging unconditionally in the space $\HH{-\tau}$, and similarly,
\[
g=\sum_{\ell=1}^\mphi \sum_{k\in \Z} \la g, \tilde \phi^\ell_{\df^J;k}\ra \phi^\ell_{\df^J;k}+ \sum_{j=J}^\infty \sum_{\ell=1}^\mpsi \sum_{k\in \Z}
\la g, \tilde \psi^\ell_{\df^j;k}\ra  \psi^\ell_{\df^j;k}, \qquad g\in \HH{\tau}
\]
with the series converging unconditionally in the space $\HH{\tau}$. See \cite{HanShen:ca:2009} for more detail.

Let $\B$ be a normed function space and $\B'$ be its dual.
For example, for Besov spaces $\BS{\tau}{p}{q}$, $(\BS{\tau}{p}{q})'=\BS{-\tau}{p'}{q'}$, where $\tau\in \R, 1\le p<\infty$ and $1\le q<\infty$ with $1/p+1/p'=1/q+1/q'=1$. Similarly, for Triebel-Lizorkin spaces $\TL{\tau}{p}{q}$, $(\TL{\tau}{p}{q})'=\TL{-\tau}{p'}{q'}$.
Replacing $\HH{\tau}$ and $\HH{-\tau}$ by $\B$ and $\B'$, respectively, the notion of a pair of nonhomogeneous dual wavelet frames can be generalized from a pair of dual Sobolev spaces $(\HH{\tau}, \HH{-\tau})$ to a pair of dual function spaces $(\B, \B')$ by a proper choice of some normed sequence spaces $\seqspace_\B$ and $\seqspace_{\B'}$. We shall not further address this issue in this paper.

In the rest of this section, we shall use the following notation
\begin{equation}\label{hatB}
\wh{\HH{\tau}}:=\{ \hat f\, : \, f\in \HH{\tau}\}\quad \hbox{and} \quad \|\hat f\|_{\wh{\HH{\tau}}}:=\|f\|_{\HH{\tau}}, \qquad f\in \HH{\tau}.
\end{equation}
%


By the following result, we see that the notion of a frequency-based nonhomogeneous dual wavelet frames in the distribution space plays a basic role in the study of pairs of nonhomogeneous dual wavelet frames in a pair of dual function spaces.

\begin{theorem}\label{thm:sobolev}
Let $\tau\in \R$ and $\df$ be a nonzero real number. Let $\Phi, \Psi, \tilde \Phi, \tilde \Psi$ in \eqref{Phi:Psi} be subsets of tempered distributions. Define
\begin{equation}\label{PhuPsi}
\fPhi:=\{\hat \phi\; : \; \phi \in \Phi\}, \quad
\fPsi:=\{\hat \psi\; : \; \psi \in \Psi\}, \quad
\tilde \fPhi:=\{\hat{\tilde \phi}\; : \; \tilde \phi \in \tilde \Phi\}, \quad
\tilde \fPsi:=\{\hat{\tilde\psi}\; : \; \tilde \psi \in \tilde \Psi\}.
\end{equation}
Then $(\WS_J(\Phi; \Psi), \WS_J(\tilde \Phi; \tilde \Psi))$ forms a pair of nonhomogeneous dual wavelet frames in the pair of dual Sobolev spaces $(\HH{\tau}, \HH{-\tau})$, if and only if,
\begin{enumerate}
\item[{(i)}] there exists a positive constant $C$ such that
\begin{equation}\label{boundedness:primal}
\big \| \{ \la \ff, \fh\ra \}_{\fh\in \FWS_J(\fPhi; \fPsi)}\big \|_{\seqspace_{\HH{\tau}}} \le C \|\ff\|_{\wh{\HH{-\tau}}}, \qquad \ff\in \D
\end{equation}
and
\begin{equation}\label{boundedness:dual}
\big \| \{ \la \ff, \tilde \fh\ra \}_{\tilde \fh\in \FWS_J(\tilde \fPhi; \tilde \fPsi)}\big \|_{\seqspace_{\HH{-\tau}}} \le C \|\ff\|_{\wh{\HH{\tau}}}, \qquad \ff\in \D;
\end{equation}

\item[{(ii)}] the pair $(\FWS_J(\fPhi; \fPsi), \FWS_J(\tilde \fPhi; \tilde \fPsi))$, which is the image of the pair $(\WS_J(\Phi; \Psi), \WS_J(\tilde \Phi; \tilde \Psi))$ under the Fourier transform, forms a pair of frequency-based nonhomogeneous dual wavelet frames in the distribution space $\Dpr$.
\end{enumerate}
\end{theorem}

\begin{proof} For $f\in \check{\mathscr{D}}(\R):=\{ h\, : \, \hat h\in \D\}$ and $g\in \HH{\alpha}$ with $\alpha\in \R$, the relation $\la f, g\ra=\frac{1}{2\pi} \la \hat f, \hat g\ra$ holds.
Since $\check{\mathscr{D}}(\R)$ is contained in both $\HH{\tau}$ and $\HH{-\tau}$ (or equivalently, $\D \subseteq \wh{\HH{\tau}}$ and $\D\subseteq \wh{\HH{-\tau}}$),
the necessity part is evident.
Hence, it suffices to prove the sufficiency part.

By \eqref{boundedness:primal}, for all $\ff\in \D$ and $N\in \N$,
\begin{equation}\label{boundedness:primal:N}
\sum_{\fphi\in \fPhi} \sum_{k=-N}^N |\df|^{-2\tau J} |\la \ff, \fphi_{\df^{-J};0,k}\ra|^2+\sum_{j=J}^N  \sum_{\fpsi\in \fPsi} \sum_{k=-N}^N |\df|^{-2\tau j}|\la \ff, \fpsi_{\df^{-j};0,k}\ra|^2\le C^2 \| \ff \|^2_{\wh{\HH{-\tau}}}.
\end{equation}
Using \eqref{boundedness:primal:N}, we now prove that $\fPhi, \fPsi\subseteq \wh{\HH{\tau}}$. By \eqref{boundedness:primal:N}, for $\fpsi\in \fPsi$, we have $|\la \ff, \fpsi_{\df^{-J};0,0}\ra|\le C |\df|^{\tau J} \|\ff\|_{\wh{\HH{-\tau}}}$ for all $\ff\in \D$. Therefore, $\la \cdot, \fpsi_{\df^{-J};0,0}\ra$ can be extended into a continuous linear functional on $\wh{\HH{-\tau}}$.
Since $\wh{\HH{\tau}}$ is the dual space of $\wh{\HH{-\tau}}$, there exists $\mathring{\fpsi}\in\wh{\HH{\tau}}$ such that $\la \cdot, \fpsi_{\df^{-J};0, 0}\ra=\la \cdot, \mathring{\fpsi}\ra$. In particular, $\la \ff, \fpsi_{\df^{-J};0,0}-\mathring{\fpsi}\ra=0$ for all $\ff\in \D$.
Since $\fpsi_{\df^{-J};0,0}-\mathring{\fpsi}$ is a distribution, we must have $\fpsi_{\df^{-J};0,0}=\mathring{\fpsi}$ in the sense of distributions. By $\mathring{\fpsi}\in \wh{\HH{\tau}}$, we deduce that $\fpsi_{\df^{-J};0,0}\in \wh{\HH{\tau}}$ and therefor, $\fpsi\in \wh{\HH{\tau}}$. Consequently, by \eqref{boundedness:primal:N}, we proved that
$\fPhi, \fPsi\subseteq \wh{\HH{\tau}}$.
In other words, we proved that
$\Phi, \Psi\subseteq \HH{\tau}$.

Since $\D$ is dense in $\wh{\HH{-\tau}}$ and since all $\la \cdot, \fphi_{\df^{-J};0,k}\ra, \fphi\in \fPhi$ and $\la \cdot, \fpsi_{\df^{-j};0,k}\ra, \fpsi\in \fPsi$ are continuous linear functionals on $\wh{\HH{-\tau}}$, we see that \eqref{boundedness:primal:N} holds for all $\ff \in \wh{\HH{-\tau}}$ and all $N\in \N$. Letting $N\to +\infty$ in \eqref{boundedness:primal:N} and noting $\la f, g\ra=\frac{1}{2\pi} \la \hat f, \hat g\ra$ for all $f\in \HH{-\tau}$ and $g\in \HH{\tau}$, we conclude that $\Phi, \Psi\subseteq \HH{\tau}$ and
\begin{equation}\label{boundedness:primal:0}
\sum_{\ell=1}^\mphi \sum_{k\in \Z} |\df|^{-2\tau J} |\la f, \phi^{\ell}_{\df^J;k}\ra|^2+\sum_{j=J}^\infty \sum_{\ell=1}^\mpsi \sum_{k\in \Z} |\df|^{-2\tau j}|\la f, \psi^\ell_{\df^j;k}\ra|^2\le \frac{C^2}{2\pi} \| f\|^2_{\HH{-\tau}}, \quad f\in \HH{-\tau}.
\end{equation}
Similarly, we can show that \eqref{boundedness:dual} implies $\tilde \Phi, \tilde \Psi \subseteq \HH{-\tau}$ and
\begin{equation}\label{boundedness:dual:0}
\sum_{\ell=1}^\mphi \sum_{k\in \Z} |\df|^{2\tau J} |\la g, \tilde \phi^{\ell}_{\df^J;k}\ra|^2+\sum_{j=J}^\infty \sum_{\ell=1}^\mpsi \sum_{k\in \Z} |\df|^{2\tau j}|\la g, \tilde \psi^\ell_{\df^j;k}\ra|^2\le \frac{C^2}{2\pi} \| g\|^2_{\HH{\tau}}, \quad g\in \HH{\tau}.
\end{equation}
Define two operators $\Wop: \HH{-\tau} \rightarrow \seqspace_{\HH{\tau}}$ and
$\tilde \Wop: \HH{\tau}\rightarrow \seqspace_{\HH{-\tau}}$ by
\begin{equation}\label{op:W}
\Wop  f:=\{ \la f, h\ra\}_{h\in \WS_J(\Phi; \Psi)} \quad \hbox{and}\quad
\tilde \Wop  g:=\{ \la g, \tilde h\ra\}_{\tilde h\in \WS_J(\tilde \Phi; \tilde \Psi)}.
\end{equation}
Then \eqref{boundedness:primal:0} and \eqref{boundedness:dual:0} are equivalent to saying that the operators $\Wop$ and $\tilde \Wop$ are well-defined bounded linear operators, more precisely,
\begin{equation}\label{bounded:W}
\|\Wop f \|_{\seqspace_{\HH{\tau}}}\le \frac{C}{\sqrt{2\pi}} \|f\|_{\HH{-\tau}}, \;\;
\|\Wop g\|_{\seqspace_{\HH{-\tau}}}\le \frac{C}{\sqrt{2\pi}} \|g\|_{\HH{\tau}} \; \forall\, f\in \HH{-\tau}, g\in \HH{\tau}.
\end{equation}
By Cauchy-Schwarz inequality, using \eqref{boundedness:primal:0} and \eqref{boundedness:dual:0}, we see that for $f\in \HH{-\tau}$ and $g\in \HH{\tau}$, the series on the right-hand side of \eqref{identity} converges absolutely. Now by assumption in item (ii), we see that \eqref{identity} holds for all $f, g\in \check{\mathscr{D}}(\R)$. In other words, for $f, g\in \check{\mathscr{D}}(\R)$,
\begin{equation}\label{s3:eq1}
\la \Wop f, \tilde \Wop g\ra
=\la f, g\ra.
\end{equation}
Note that $\check{\mathscr{D}}(\R)$ is dense in both $\HH{-\tau}$ and $\HH{\tau}$. Now we use a standard argument to show that \eqref{s3:eq1} holds for all $f\in \HH{-\tau}$ and $g\in \HH{\tau}$. For $f\in \HH{-\tau}$ and $g\in \HH{\tau}$, there exist two sequences $\{ f_n\}_{n=1}^\infty$ and $\{g_n\}_{n=1}^\infty$ in $\check{\mathscr{D}}(\R)$ such that
\begin{equation}\label{approx:f:g}
\lim_{n\to \infty} \| f_n-f\|_{\HH{-\tau}}=0 \quad \hbox{and}\quad \lim_{n\to \infty} \|g_n -g\|_{\HH{\tau}}=0.
\end{equation}
Observe that
\[
\la \Wop f, \tilde \Wop g\ra=\la \Wop(f-f_n), \tilde \Wop g\ra+\la \Wop  f_n, \tilde \Wop (g-g_n)\ra+\la \Wop f_n, \tilde \Wop g_n\ra.
\]
By $f_n, g_n\in \check{\mathscr{D}}(\R)$ and \eqref{s3:eq1}, we have $\la \Wop f_n, \tilde \Wop g_n\ra=\la f_n, g_n\ra$. Therefore, we have
\[
\la \Wop f, \tilde \Wop g\ra-\la f, g\ra=\la \Wop(f-f_n), \tilde \Wop g\ra+\la \Wop f_n, \tilde \Wop (g-g_n)\ra+\la f_n-f, g_n\ra+\la f, g_n-g\ra.
\]
Now by \eqref{bounded:W} and the triangle inequality,
\begin{align*}
|\la \Wop f, \tilde \Wop g\ra-\la f, g\ra|
&\le |\la \Wop(f-f_n), \tilde \Wop g\ra|+|\la \Wop f_n, \tilde \Wop (g-g_n)\ra|+|\la f_n-f, g_n\ra|+|\la f, g_n-g\ra|\\
&\le C_1\Big( \|f-f_n\|_{\HH{-\tau}}\|g\|_{\HH{\tau}}
+\|f_n\|_{\HH{-\tau}}\|g-g_n\|_{\HH{\tau}}\\
&\quad +\|f_n-f\|_{\HH{-\tau}}\|g_n\|_{\HH{\tau}}
+\|f\|_{\HH{-\tau}}\|g_n-g\|_{\HH{\tau}}\Big),
\end{align*}
where $C_1:=\frac{C^2}{2\pi}+1$. By \eqref{approx:f:g}, the right-hand side of the above inequality goes to $0$ as $n\to \infty$. Consequently, \eqref{s3:eq1} holds for all $f\in \HH{-\tau}$ and $g\in \HH{\tau}$. That is, \eqref{identity} and item (3) have been verified.

To prove item (1), by \eqref{s3:eq1} and \eqref{bounded:W}, we have
\[
|\la f, g\ra|=|\la \Wop f, \tilde \Wop g\ra|\le \| \Wop f\|_{\seqspace_{\HH{\tau}}} \|\tilde \Wop g\|_{\seqspace_{\HH{-\tau}}} \le \frac{C}{\sqrt{2\pi}} \|g\|_{\HH{\tau}}  \| \Wop f\|_{\seqspace_{\HH{\tau}}}.
\]
That is, for $f\in \HH{-\tau}$, the following inequality holds:
\[
\| \Wop f\|_{\seqspace_{\HH{\tau}}} \ge \frac{\sqrt{2\pi}}{C} \sup_{g\in \HH{\tau} \bs \{0\}} \frac{|\la f, g\ra|}{\|g\|_{\HH{\tau} }}=\frac{\sqrt{2\pi}}{C} \|f\|_{\HH{\tau}},
\]
where in the last step we used the fact that $\HH{\tau}$ is the dual space of $\HH{-\tau}$. Hence, item (1) holds. Item (2) can be proved similarly. This completes the proof of the sufficiency part.
\end{proof}

Next, we discuss connections between nonhomogeneous and homogeneous  wavelet systems in the particular function space $\Lp{2}$.
To do so, we need an auxiliary result, which is essentially known in the literature,  e.g., see \cite[Page~28]{Han:oep}. For completeness, we present a proof here.

\begin{lemma}\label{lem:to0} Let $\df$ be a real number such that $|\df|>1$.
Let $\phi\in \Lp{2}$ such that there exists a positive constant $C$ such that $\sum_{k\in \Z} |\la f, \phi(\cdot-k)\ra|^2\le C\|f\|_{\Lp{2}}^2$ for all $f\in \Lp{2}$ (or equivalently, $\sum_{k\in \Z} |\hat\phi(\xi+2\pi k)|^2\le C$ for almost every $\xi\in \R$). Then
\begin{equation}\label{lem:to0:eq1}
\lim_{j\to -\infty} \sum_{k\in \Z} |\la f, \phi_{\df^j;k}\ra|^2=0 \qquad \forall\; f\in \Lp{2}.
\end{equation}
\end{lemma}

\begin{proof} It suffices to prove the case $\df>1$, since the negative case $\df<-1$ can be proved similarly. So, we assume $\df>1$.
We first show that \eqref{lem:to0:eq1} holds for all $f=\chi_{(t_1, t_2)}$ with $t_1<t_2$, where $\chi_{(t_1, t_2)}$ denotes the characteristic function of the open interval $(t_1, t_2)$. For $f=\chi_{(t_1, t_2)}$, by Cauchy-Schwarz inequality, we have
\[
|\la f, \phi_{\df^j;k}\ra|^2
=\df^{-j} \left| \int_{(\df^j t_1-k, \df^j t_2-k)} \phi(x) dx\right|^2\le
\df^{-j} \Big(\int_{(\df^j t_1-k, \df^j t_2-k)} dx\Big)
\Big( \int_{(\df^j t_1-k, \df^j t_2-k)} |\phi(x)|^2 dx\Big).
\]
Since $\int_{(\df^j t_1-k, \df^j t_2-k)} dx=\df^j (t_2-t_1)$, noting that $\lim_{j\to -\infty}\df^j=0$ by $\df>1$ and all $(\df^j t_1-k, \df^j t_2-k), k\in \Z$ are disjoint as $j\to -\infty$, we deduce that
\[
\sum_{k\in \Z} |\la f, \phi_{\df^j;k}\ra|^2\le |t_2-t_1| \int_{\cup_{k\in \Z} (\df^j t_1-k, \df^j t_2-k)} |\phi(x)|^2 dx \to 0,
\]
as $j\to -\infty$, since $\phi\in \Lp{2}$.
Consequently, \eqref{lem:to0:eq1} holds for all $f=\chi_{(t_1, t_2)}$. So, \eqref{lem:to0:eq1} holds for all $f$ that are finite linear combinations of characteristic functions of bounded open intervals.

On the other hand, define operators $P_j: \Lp{2} \rightarrow \lp{2}$ by $P_j f:=\{ \la f, \phi_{\df^j;k}\ra\}_{k\in \Z}$. By $P_j f=P_0 f_{\df^{-j};0,0}$, we have
\[
\|P_j f\|_{\lp{2}}^2=\|P_0 f_{\df^{-j};0,0}\|_{\lp{2}}^2\le C \| f_{\df^{-j};0,0}\|_{\Lp{2}}^2=
C\|f\|_{\Lp{2}}^2.
\]
When $\sum_{k\in \Z} |\hat\phi(\xi+2\pi k)|^2\le C$ for almost every $\xi\in \R$, by Lemma~\ref{lem:converg}, we see that
\begin{align*}
\|P_j f\|^2_{\lp{2}}&=\frac{1}{4\pi^2}\sum_{k\in \Z}|\la \hat{f}, \wh{\phi_{\df^j;k}}\ra|^2
=\frac{\df^j}{2\pi}
\int_{-\pi}^\pi \Big| \sum_{k\in \Z} \hat f(\df^j(\xi+2\pi k))\ol{\hat \phi(\xi+2\pi k)}\Big|^2 d\xi\\
&\le\frac{\df^j}{2\pi} \int_{-\pi}^\pi \Big( \sum_{k\in \Z} |\hat f(\df^j(\xi+2\pi k))|^2\Big)\Big(\sum_{k\in \Z} |\hat \phi(\xi+2\pi k)|^2\Big) d\xi\\
&\le \frac{C}{2\pi} \df^j \int_{-\pi}^\pi \sum_{k\in \Z} |\hat f (\df^j(\xi+2\pi k))|^2 d \xi=\frac{C}{2\pi} \int_\R |\hat f(\xi)|^2=C \|f\|^2_{\Lp{2}}.
\end{align*}
%

Let $f\in \Lp{2}$. For an arbitrary $\gep>0$, there is $g$, which is a finite linear combination of characteristic functions of bounded open intervals, such that $\|f-g\|_{\Lp{2}}\le \gep$. Since $\lim_{j\to -\infty} \| P_j g\|_{\lp{2}}=0$, there exists $J$ such that $\|P_j g\|_{\lp{2}}\le \gep$ for all $j\le J$. Hence, for all $j\le J$,
we have
\[
\|P_j f\|_{\lp{2}}\le \| P_j (f-g)\|_{\lp{2}}+\| P_j g\|_{\lp{2}}\le \sqrt{C} \| f-g\|_{\Lp{2}}+
\| P_j g\|_{\lp{2}}\le (\sqrt{C}+1) \gep.
\]
Hence, $\lim_{j\to -\infty} \|P_j f\|_{\lp{2}}=0$. That is, \eqref{lem:to0:eq1} holds.
\end{proof}

\begin{theorem}\label{thm:connection}
Let $\df$ be a real number. Let $\Phi, \Psi, \tilde \Phi, \tilde \Psi$ in \eqref{Phi:Psi} be subsets of $\Lp{2}$.
Suppose that $\WS_J(\Phi; \Psi)$, which is defined in \eqref{WS}, is a frame in $\Lp{2}$ for some integer $J$, that is, there are positive constants $C_1$ and $C_2$ such that
\begin{equation}\label{frameinL2}
C_1 \|f\|^2_{\Lp{2}}\le  \sum_{\ell=1}^\mphi \sum_{k\in \Z} |\la f, \phi^\ell_{\df^J;k}\ra|^2+\sum_{j=J}^\infty  \sum_{\ell=1}^\mpsi \sum_{k\in \Z} |\la f, \psi^\ell_{\df^j;k}\ra|^2\le C_2 \|f\|^2_{\Lp{2}}, \qquad f\in \Lp{2}.
\end{equation}
Then \eqref{frameinL2} holds for all integers $J$.
If in addition $|\df|>1$, then $\WS(\Psi)$ must be a frame in $\Lp{2}$ with the same frame bounds satisfying
\begin{equation}\label{frameinL2:hom}
C_1 \|f\|^2_{\Lp{2}}\le \sum_{j\in \Z} \sum_{\ell=1}^\mpsi \sum_{k\in \Z} |\la f, \psi^\ell_{\df^j;k}\ra|^2\le C_2 \|f\|^2_{\Lp{2}}, \qquad f\in \Lp{2}.
\end{equation}
If $|\df|>1$ and the pair $(\WS_J(\Phi; \Psi), \WS_J(\tilde \Phi, \tilde \Psi))$ forms a pair of nonhomogeneous dual wavelet frames in $(\Lp{2}, \Lp{2})$ for some integer $J$, then $(\WS(\Psi), \WS(\tilde \Psi))$ forms a pair of homogeneous dual wavelet frames in $\Lp{2}$, that is,
both $\WS(\Psi)$ and $\WS(\tilde \Psi)$ are frames in $\Lp{2}$, and the following identity holds
\begin{equation}\label{hom:dwf}
\la f, g\ra=\sum_{j\in \Z}  \sum_{\ell=1}^\mpsi \sum_{k\in \Z} \la f, \psi^\ell_{\df^j;k}\ra \la \tilde \psi^\ell_{\df^j;k}, g\ra, \qquad f,g\in \Lp{2}
\end{equation}
with the series converging absolutely.
\end{theorem}

\begin{proof} By the simple observation in \eqref{rel:dilation}, it is easy to see that for all integers $J$, \eqref{frameinL2} holds with the same constants $C_1$ and $C_2$.
In particular, for all $J\in \Z$, we have
%
\[
\|P_J f\|^2:=\sum_{\ell=1}^\mphi \sum_{k\in \Z} |\la f, \phi^\ell_{\df^J;k}\ra|^2\le C_2 \|f\|^2_{\Lp{2}}.
\]
%
Since $|\df|>1$, by Lemma~\ref{lem:to0},  we have $\lim_{J\to -\infty}\sum_{\ell=1}^\mphi\sum_{k\in \Z} |\la f, \phi^\ell_{\df^J;k}\ra|^2=0$ for all $f\in \Lp{2}$. Now it is easy to deduce that \eqref{frameinL2:hom} holds.

To prove the second claim, by what has been proved,
both $\WS(\Psi)$ and
$\WS(\tilde \Psi)$ are frames in $\Lp{2}$.
Note that the pair $(\WS_J(\Phi; \Psi), \WS_J(\tilde \Phi, \tilde \Psi))$ forms a pair of nonhomogeneous dual wavelet frames in $(\Lp{2}, \Lp{2})$ for all integers $J$. Thus, for all integers $J$, we have
\begin{equation}\label{dwf:eq1}
\la f, g\ra=\sum_{\ell=1}^\mphi\sum_{k\in \Z} \la f, \phi^\ell_{\df^J;k}\ra\la \tilde \phi^\ell_{\df^J;k}, g\ra+
\sum_{j=J}^\infty  \sum_{\ell=1}^\mpsi \sum_{k\in \Z} \la f, \psi^\ell_{\df^j;k}\ra \la \tilde \psi^\ell_{\df^j;k}, g\ra, \qquad f,g\in \Lp{2}
\end{equation}
with the series on the right-hand side converging absolutely.
By Lemma~\ref{lem:to0} again, for $f,g\in \Lp{2}$, we have $\lim_{J\to -\infty}\sum_{\ell=1}^\mphi\sum_{k\in \Z} \la f, \phi^\ell_{\df^J;k}\ra\la \tilde \phi^\ell_{\df^J;k}, g\ra=0$. Now we see that \eqref{hom:dwf} holds.
\end{proof}

As a direct consequence of all the above results, we have

\begin{cor}\label{cor:twf} Let $\df$ be a real number such that $|\df|>1$. Let $\fPhi=\{\fphi^1, \ldots, \fphi^\mphi\}$ and $\fPsi=\{\fpsi^1, \ldots, \fpsi^\mpsi\}$ be subsets of distributions in $\Dpr$. Then the following statements are equivalent
\begin{enumerate}
\item $\FWS_J(\fPhi; \fPsi)$, which is defined in \eqref{FWS}, is a frequency-based nonhomogeneous tight wavelet frame in $\Lp{2}$ for some integer $J$, that is, $\fPhi, \fPsi\subseteq \Lp{2}$ and
\begin{equation}\label{twf:2}
\sum_{\ell=1}^\mphi \sum_{k\in \Z} |\la \ff, \fphi^{\ell}_{\df^{-J};0,k} \ra|^2+\sum_{j=J}^\infty\sum_{\ell=1}^\mpsi  \sum_{k\in \Z} | \la \ff, \fpsi^{\ell}_{\df^{-j};0,k} \ra|^2=2\pi \|\ff\|^2_{\Lp{2}}\qquad \forall\; \ff\in \Lp{2};
\end{equation}

\item item (1) is true for all integers $J$;

\item $(\FWS_J(\fPhi; \fPsi), \FWS_J(\fPhi; \fPsi))$ forms a pair of frequency-based nonhomogeneous dual wavelet frames in the distribution space $\Dpr$;
\item $\fPhi, \fPsi\subseteq \lLp{2}$ and \eqref{main:eq1:dm}, \eqref{I:eq1:sp}, \eqref{I:eq2:sp}, \eqref{I:eq3:sp} hold with $\tilde \fPhi:=\fPhi$ and $\tilde \fPsi:=\fPsi$;
\item there exist $\phi^1, \ldots, \phi^\mphi, \psi^1, \ldots, \psi^\mpsi\in \Lp{2}$ such that $\fphi^1=\wh{\phi^1}, \ldots, \fphi^\mphi=\wh{\phi^\mphi}, \fpsi^1=\wh{\psi^1}, \ldots, \fpsi^\mpsi=\wh{\psi^\mpsi}$, and $\WS_J(\{\phi^1, \ldots, \phi^\mphi\}; \{ \psi^1, \ldots, \psi^\mpsi\})$ is a nonhomogeneous tight wavelet frame in $\Lp{2}$:
\begin{equation}\label{twf}
\sum_{\ell=1}^\mphi \sum_{k\in \Z} |\la f, \phi^{\ell}_{\df^J;k} \ra|^2+\sum_{j=J}^\infty\sum_{\ell=1}^\mpsi  \sum_{k\in \Z} | \la f, \psi^{\ell}_{\df^j;k} \ra|^2=\|f\|^2_{\Lp{2}}\qquad \forall\; f\in \Lp{2}.
\end{equation}
\end{enumerate}
Moreover, any of the above statements implies that $\FWS(\fPsi)$ is a frequency-based homogeneous tight wavelet frame in $\Lp{2}$, that is, $\fPsi\subseteq \Lp{2}$ and
\begin{equation}\label{htwf:2}
\sum_{j\in \Z}  \sum_{\ell=1}^\mpsi \sum_{k\in \Z} | \la \ff, \fpsi^{\ell}_{\df^{j};0,k}\ra|^2= 2\pi \|\ff\|^2_{\Lp{2}}\qquad \forall\; \ff\in \Lp{2}.
\end{equation}
\end{cor}

\begin{proof} We first show that if item (3) holds, then
\eqref{boundedness:primal} must be true with $C=2\pi$ and $\tau=0$. Since item (3) holds, by definition, for all $\ff\in \D$, we have
\[
\lim_{J'\to +\infty}
\Big( \sum_{\ell=1}^\mphi \sum_{k\in \Z} |\la \ff, \fphi^\ell_{\df^{-J}; 0, k}\ra|^2+\sum_{j=J}^{J'-1}\sum_{\ell=1}^\mpsi \sum_{k\in \Z}
|\la \ff, \fpsi^\ell_{\df^{-j}; 0, k}\ra|^2\Big)=2\pi \|\ff\|_{\Lp{2}}^2.
\]
Now it is straightforward to see that \eqref{boundedness:primal} must be true with $C=2\pi$ and $\tau=0$, since $\|\f\|_{\wh{\Lp{2}}}^2=\frac{1}{2\pi} \|\ff\|_{\Lp{2}}^2$. Consequently, by Theorem~\ref{thm:sobolev}, item (3) implies $\fPhi, \fPsi\subseteq \Lp{2}$.

Now by Theorem~\ref{thm:sobolev}, (1), (2), and (3) are equivalent to each other. The equivalence of (3) and (4) is guaranteed by Theorem~\ref{thm:main:1:general}. The equivalence between (1) and (5) is trivial. \eqref{htwf:2} is a direct consequence of Theorem~\ref{thm:connection}.
\end{proof}

From the following result, we see that there is a natural connection between refinable function vectors and
frequency-based nonhomogeneous orthonormal wavelet bases in $\Lp{2}$.

\begin{prop}\label{prop:orth} Let $\df$ be a nonzero real number. Let $\fPhi$ and $\fPsi$ in \eqref{gen:set} be subsets of $\Lp{2}$. Suppose that
$\frac{1}{\sqrt{2\pi}} \FWS_J(\fPhi; \fPsi)$ is an orthonormal basis of $\Lp{2}$ for some integer $J$ (this is equivalent to saying that \eqref{twf:2} holds with $\|\fphi^1\|_{\Lp{2}}=\cdots=\|\fphi^\mphi\|_{\Lp{2}}
=\|\fpsi^1\|_{\Lp{2}}=\cdots=\|\fpsi^\mpsi\|_{\Lp{2}}=
\sqrt{2\pi}$).
Denote $\vec{\fphi}:=(\fphi^1, \ldots, \fphi^\mphi)^T$ and $\vec{\fpsi}:=(\fpsi^1, \ldots, \fpsi^\mpsi)^T$. Then there must exist $\mphi\times \mphi$ matrix $\fa$ and $\mpsi \times \mphi$ matrix $\fb$ of $2\pi$-periodic measurable functions in $\TLp{2}$ such that
\begin{equation}\label{rel:refeq}
\vec{\fphi}(\df\xi)=\fa(\xi)\vec{\fphi}(\xi) \quad \hbox{and}\quad
\vec{\fpsi}(\df\xi)=\fb(\xi) \vec{\fphi}(\xi), \qquad a.e.\; \xi\in \R.
\end{equation}
Moreover, if $|\df|>1$, then $\frac{1}{\sqrt{2\pi}} \FWS(\fPsi)$ is also an orthonormal basis of $\Lp{2}$.
\end{prop}

\begin{proof} Note that $\frac{1}{\sqrt{2\pi}} \FWS_J(\fPhi; \fPsi)$ is an orthonormal basis of $\Lp{2}$ if and only if it is an orthonormal basis of $\Lp{2}$ for all integers $J$. Consider the expansion of the elements $\fphi^{\ell'}_{\df;0,0}$ and $\fpsi^{\ell'}_{\df;0,0}$ under the orthonormal basis $\frac{1}{\sqrt{2\pi}} \FWS_0(\fPhi; \fPsi)$. By orthogonality, we have $\fphi^{\ell'}_{\df;0,0}=\frac{1}{2\pi}\sum_{\ell=1}^\mphi \sum_{k\in \Z} \la \fphi^{\ell'}_{\df;0,0}, \fphi^{\ell}_{1;0,k}\ra \fphi^{\ell}_{1; 0, k}$ and $\fpsi^{\ell'}_{\df;0,0}=\frac{1}{2\pi}\sum_{\ell=1}^\mphi \sum_{k\in \Z} \la \fpsi^{\ell'}_{\df;0,0}, \fphi^{\ell}_{1;0,k}\ra \fphi^{\ell}_{1; 0, k}$.
Noting that $\fphi^{\ell'}_{\df;0,0}(\xi)=|\df|^{1/2} \fphi^{\ell'}(\df \xi)$ and $\fphi^{\ell}_{1;0,k}(\xi)=e^{-ik\xi} \fphi^\ell(\xi)$,
we see that \eqref{rel:refeq} holds with
\[
[\fa(\xi)]_{\ell',\ell}=\frac{1}{2\pi\sqrt{|\df|}}\sum_{k\in \Z} \la \fphi^{\ell'}_{\df;0,0}, \fphi^{\ell}_{1;0,k}\ra e^{-i k\xi}, \qquad \ell, \ell'=1, \ldots, \mphi
\]
and
\[
[\fb(\xi)]_{\ell',\ell}=\frac{1}{2\pi\sqrt{|\df|}}\sum_{k\in \Z} \la \fpsi^{\ell'}_{\df;0,0}, \fphi^{\ell}_{1;0,k}\ra e^{-i k\xi}, \qquad \ell=1, \ldots, \mphi,\, \ell'=1, \ldots, \mpsi,
\]
where $[\fa(\xi)]_{\ell', \ell}$ denotes the $(\ell',\ell)$-entry of the matrix $\fa(\xi)$. Since $\sum_{k\in \Z} |\la \fphi^{\ell'}_{\df;0,0}, \fphi^{\ell}_{1;0, k}\ra|^2<\infty$,
all $[\fa]_{\ell',\ell}$ are well defined elements in $\TLp{2}$.
Similarly, all $[\fb]_{\ell',\ell}$ are well defined elements in $\TLp{2}$.
\end{proof}


\section{Nonstationary Dual Wavelet Frames in the Distribution Space}

All the results in the previous sections have been mainly built on the multiresolution-like structure in \eqref{MRA} for stationary nonhomogeneous wavelet systems. Here stationary means that at scale level $j$ the dilation is $\df^j$ and the generating wavelet functions are independent of the scale level $j$. The result in Lemma~\ref{lem:twolevel} characterizing the multiresolution-like structure in \eqref{twolevel:general} makes most proofs  in the previous sections relatively simple. Nonhomogeneous wavelet systems are closely related to nonstationary wavelets, which are useful in many applications since the nonstationary wavelet filter banks can be implemented in almost the same way and efficiency as a traditional fast wavelet transform. However, except a few special cases as discussed in \cite{CHS:ntwf,CohenDyn,HanShen:siam:2008} and some references therein, only few theoretical results on nonstationary wavelets are available in the literature.

In this section, we shall see that the notion of a pair of frequency-based nonhomogeneous dual wavelet frames in the distribution space is very flexible and similar results hold in the most general setting of fully nonstationary wavelets. Since there are few theoretical results on nonstationary wavelets in the literature,
it is worth our effort to provide a better picture to understand them in this section.

Let us first introduce the notion of a pair of frequency-based nonstationary dual wavelet frames in the distribution space.
Let $J\in \Z$ and $\{\gl_j\}_{j=J}^\infty$ be a sequence of nonzero real numbers. Let
\begin{equation}\label{ns:Phi}
\fPhi=\{\fphi^1, \ldots, \fphi^\mphi\}, \quad
\tilde \fPhi=\{\tilde \fphi^1, \ldots, \tilde \fphi^\mphi\}
\end{equation}
and
\begin{equation}\label{ns:Psi}
\fPsi^j:=\{\fpsi^{j,1}, \ldots, \fpsi^{j, \mpsi_j}\},\qquad
\tilde \fPsi^j:=\{\tilde \fpsi^{j,1}, \ldots, \tilde \fpsi^{j, \mpsi_j}\}
\end{equation}
be subsets of distributions in $\Dpr$ with $j\ge J$ and $\mpsi_j\in \N$. We say that the pair
\begin{equation}\label{FWS:ns}
(\FWS_J(\fPhi; \{ \fPsi^{j}\}_{j=J}^{\infty}), \FWS_J(\tilde \fPhi; \{ \tilde \fPsi^{j}\}_{j=J}^{\infty}))
\end{equation}
forms {\it a pair of frequency-based nonstationary dual wavelet frames in the distribution space $\Dpr$} if the following identity holds:
\begin{equation}\label{ndwf}
\sum_{\ell=1}^{\mphi} \sum_{k\in \Z} \la \ff, \fphi^\ell_{\gl_J;0,k}\ra \la \tilde \fphi^\ell_{\gl_J;0,k}, \fg\ra+
 \sum_{j=J}^\infty  \sum_{\ell=1}^{\mpsi_j} \sum_{k\in \Z}
\la \ff, \fpsi^{j,\ell}_{\gl_j;0,k}\ra \la \tilde \fpsi^{j,\ell}_{\gl_j; 0,k}, \fg\ra=2\pi \la \ff, \fg\ra \quad \forall\; \ff, \fg\in \D,
\end{equation}
where the infinite series in \eqref{ndwf} converge in the following sense:
\begin{enumerate}
\item For every $\ff,\fg\in \D$, all the series $\sum_{k\in \Z} \la \ff, \fphi^\ell_{\gl_J; 0,k}\ra \la \tilde \fphi^\ell_{\gl_J; 0,k}, \fg\ra$ and $\sum_{k\in \Z}
\la \ff, \fpsi^{j,\ell'}_{\gl_j; 0,k}\ra \la \tilde \fpsi^{j,\ell'}_{\gl_j; 0,k}, \fg\ra$ converge absolutely for every integer $j\ge J$, $\ell=1, \ldots, \mphi$, and $\ell'=1, \ldots, \mpsi_j$;

\item for every $\ff, \fg\in \D$, the following limit exists and
\begin{equation}\label{def:ndwf}
\lim_{J'\to +\infty} \left(\sum_{\ell=1}^{\mphi} \sum_{k\in \Z} \la \ff, \fphi^\ell_{\gl_J; 0,k}\ra \la \tilde \fphi^{\ell}_{\gl_J; 0,k}, \fg\ra+
 \sum_{j=J}^{J'-1} \sum_{\ell=1}^{\mpsi_j}\sum_{k\in \Z}
\la \ff, \fpsi^{j,\ell}_{\gl_j; 0,k}\ra \la \tilde \fpsi^{j,\ell}_{\gl_j; 0,k}, \fg\ra\right)=2\pi \la \ff, \fg\ra.
\end{equation}
\end{enumerate}

The stationary nonhomogeneous wavelet systems considered in previous sections correspond to the case that $\gl_j=\df^{-j}$, $\mpsi_j=\mpsi$, and $\fPsi^{j}=\fPsi$ for all $j\ge J$; that is, the generating wavelet functions remain stationary (unchanged) at all the scale levels $j$.

A pair of frequency-based nonstationary dual wavelet frames in the distribution space $\Dpr$ can be similarly characterized by the following result.

\begin{theorem}\label{thm:ndwf}
Let $J$ be an integer and $\{\gl_j\}_{j=J}^\infty$ be a sequence of nonzero real numbers such that $\lim_{j\to +\infty} \gl_j=0$. Let $\fPhi, \tilde \fPhi$ in \eqref{ns:Phi} and $\fPsi^{j},\tilde \fPsi^{j}$ in \eqref{ns:Psi} be subsets of $\lLp{2}$ for all integers $j\ge J$. Then
the pair in \eqref{FWS:ns}
forms a pair of frequency-based nonstationary dual wavelet frames in the distribution space $\Dpr$, if and only if,
\begin{equation}\label{I:ndwf:eq1}
I_\fPhi^{\gl_J k}(\gl_J \xi)+\sum_{j=J}^\infty I_{\fPsi^j}^{\gl_j k}(\gl_j \xi)=0, \qquad a.e.\; \xi\in \R,\; k\in \gL\bs \{0\}
\end{equation}
(All the above infinite sums are in fact finite, since $\lim_{j\to +\infty} \gl_j k=0$ for all $k\in \R$.) and
\begin{equation}\label{I:ndwf:eq2}
\lim_{J'\to +\infty}\Big(I_\fPhi^0(\gl_J\cdot)+\sum_{j=J}^{J'-1} I_{\fPsi^j}^0(\gl_j\cdot)\Big)=1 \qquad \hbox{in the sense of distributions},
\end{equation}
where $\gL:=\cup_{j=J}^\infty [\gl_j^{-1} \Z]$ and
\begin{equation}\label{I:Phi:0}
I_\fPhi^k(\xi):=\sum_{\ell=1}^\mphi \ol{\fphi^\ell(\xi)}\tilde \fphi^\ell(\xi+2\pi k), \qquad k\in \Z \quad \hbox{and}\quad
I_\fPhi^k(\xi)\equiv 0, \qquad k\in \R \bs \Z
\end{equation}
and
\begin{equation}\label{I:Psi}
I_{\fPsi^j}^k(\xi):=\sum_{\ell=1}^{\mpsi_j} \ol{\fpsi^{j,\ell}(\xi)}\tilde \fpsi^{j,\ell}(\xi+2\pi k), \qquad k\in \Z \quad \hbox{and}\quad
I_{\fPsi^j}^k(\xi)\equiv 0, \qquad k\in \R \bs \Z.
\end{equation}
\end{theorem}

\begin{proof} Let $\ff, \fg\in \D$. By Lemma~\ref{lem:converg}, we have
\begin{equation}\label{s4:eq1}
\sum_{\ell=1}^{\mpsi_j} \sum_{k\in \Z} \la \ff, \fpsi^{j,\ell}_{\gl_j;0,k}\ra \la \tilde \fpsi^{j, \ell}_{\gl_j;0, k}, \fg\ra=
2\pi \int_\R \sum_{k\in [\gl_j^{-1} \Z]} \ff(\xi) \ol{ \fg(\xi+2\pi k)} I_{\fPsi^j}^{\gl_j k}(\gl_j\xi)\, d\xi.
\end{equation}
Since $\ff$ and $\fg$ have compact support and $\lim_{j\to +\infty} \gl_j=0$, we observe that there exists an integer $J_{\ff,\fg}$ such that $\ff(\xi) \ol{\fg(\xi+2\pi k)}\equiv 0$ for all $k\in [\gl_j^{-1} \Z]\bs \{0\}$ and $j\ge J_{\ff,\fg}$. Therefore,
\begin{equation}\label{s4:eq2}
\sum_{\ell=1}^{\mpsi_j} \sum_{k\in \Z} \la \ff, \fpsi^{j,\ell}_{\gl_j;0,k} \ra \la \tilde \fpsi^{j, \ell}_{\gl_j;0,k}, \fg\ra=
2\pi \int_\R \ff(\xi) \ol{ \fg(\xi)} I_{\fPsi^j}^{0}(\gl_j\xi)\, d\xi \qquad \forall\; j\ge J_{\ff,\fg}.
\end{equation}
Sufficiency. For $J'>J$, define
\[
S_J^{J'}(\ff,\fg):=
\sum_{\ell=1}^{\mphi} \sum_{k\in \Z} \la \ff, \fphi^{\ell}_{\gl_J;0,k}\ra \la \tilde \fphi^\ell_{\gl_J; 0,k}, \fg\ra+
 \sum_{j=J}^{J'-1} \sum_{\ell=1}^{\mpsi_j}\sum_{k\in \Z}
\la \ff, \fpsi^{j,\ell}_{\gl_j;0,k}\ra \la \tilde \fpsi^{j,\ell}_{\gl_j; 0,k}, \fg\ra.
\]
Therefore, by \eqref{s4:eq1}, for $J'>J$, we have
\begin{equation}\label{s4:eq3}
S_J^{J'}(\ff,\fg)= 2\pi \int_\R \sum_{k\in \gL} \ff(\xi) \ol{\fg(\xi+2\pi k)} \Big[ I_\fPhi^{\gl_J k}(\gl_J\xi)+\sum_{j=J}^{J'-1} I_{\fPsi^j}^{\gl_j k}(\gl_j \xi)\Big]\, d\xi.
\end{equation}
Now by \eqref{I:ndwf:eq1}, for all $J'>\max(J, J_{\ff,\fg})$, we deduce that
\begin{equation}\label{s4:eq0}
S_J^{J'}(\ff,\fg)=2\pi\int_\R \ff(\xi) \ol{\fg(\xi)} \Big[ I_\fPhi^{0}(\gl_J\xi)+\sum_{j=J}^{J'-1} I_{\fPsi^j}^{0}(\gl_j \xi)\Big]\, d\xi.
\end{equation}
Now by \eqref{I:ndwf:eq2} and \eqref{s4:eq0}, we conclude that $\lim_{J'\to +\infty} S_J^{J'}(\ff,\fg)=2\pi \int_\R \ff(\xi) \ol{\fg(\xi)} \, d\xi=2\pi \la \ff, \fg\ra$.

Necessity. The proof of the necessity part is essentially the same as that of Lemma~\ref{lem:twolevel}. Since $\lim_{j\to +\infty} \gl_j=0$, the set $\gL$ is discrete and closed.
For any temporarily fixed $\xi\in \R$ and $k_0\in \gL \bs \{0\}$, it is important to notice that $\hbox{dist}(k_0, \gL\bs \{k_0\})>0$. Now the same argument as in the proof of Lemma~\ref{lem:twolevel} leads to \eqref{I:ndwf:eq1}. Similarly, for any temporarily fixed $\xi\in \R$, since $\hbox{dist}(0, \gL\bs \{0\})>0$, by \eqref{s4:eq0}, the same argument as in the proof of Lemma~\ref{lem:twolevel} leads to \eqref{I:ndwf:eq2}.
\end{proof}

For the particular case $\gl_j=\df^{-j}$, the following result is a direct consequence of Theorem~\ref{thm:ndwf}.

\begin{cor}\label{cor:ndwf}
Let $\df$ be an integer with $|\df|>1$ and $J$ be an integer.
Let $\fPhi, \tilde \fPhi$ in \eqref{ns:Phi} and $\fPsi^j,\tilde \fPsi^{j}$ in \eqref{ns:Psi} be subsets of $\lLp{2}$ for all integers $j\ge J$. Then the pair in \eqref{FWS:ns},
with $\gl_j=\df^{-j}$ for all $j\ge J$,
forms a pair of frequency-based nonstationary dual wavelet frames in the distribution space $\Dpr$, if and only if,
for all $j_0\in \N \cup \{0\}$ and $k_0 \in \Z \bs [\df \Z]$,
\begin{equation}\label{cor:ndwf:eq1}
\sum_{\ell=1}^{\mphi}
\ol{\fphi^{\ell}(\df^{j_0}\xi)} \tilde \fphi^{\ell}(\df^{j_0}(\xi+2\pi k_0))+\sum_{j=0}^{j_0} \sum_{\ell=1}^{\mpsi_{j+J}} \ol{\fpsi^{j+J,\ell}(\df^{j_0-j}\xi)} \tilde \fpsi^{j+J,\ell}(\df^{j_0-j}(\xi+2\pi k_0))=0,
\end{equation}
for almost every $\xi\in \R$, and the following limit holds in the sense of distributions:
\begin{equation}\label{cor:ndwf:eq2}
\lim_{J'\to +\infty}\Big(\sum_{\ell=1}^{\mphi} \ol{\fphi^{\ell}(\df^{-J}\cdot)} \tilde \fphi^{\ell}(\df^{-J}\cdot)+
\sum_{j=J}^{J'-1} \sum_{\ell=1}^{\mpsi_j}  \ol{\fpsi^{j,\ell}(\df^{-j}\cdot)} \tilde \fpsi^{j,\ell}(\df^{-j}\cdot)\Big)=1.
\end{equation}
\end{cor}

\begin{proof} Since $\gl_j=\df^{-j}$, it is evident that \eqref{I:ndwf:eq2} is just \eqref{cor:ndwf:eq2}.
Note that $\gL=\cup_{j=J}^\infty [\df^j \Z]=\df^J \Z$ by $\df\in \Z$.
For any $k\in \gL\bs \{0\}=[\df^J\Z] \bs \{0\}$, we can write
$k=\df^{J+j_0} k_0$ for a unique integer $j_0\ge 0$ and a unique
$k_0\in \Z \bs [\df \Z]$. Now it is easy to check that \eqref{I:ndwf:eq1} is equivalent to \eqref{cor:ndwf:eq1}.
\end{proof}

As another application of Theorem~\ref{thm:ndwf}, using a similar argument as in the proof of Corollary~\ref{cor:twf}, we have the following result on frequency-based nonstationary tight wavelet frames in $\Lp{2}$.

\begin{cor}\label{cor:ntwf}
Let $J$ be an integer and $\{\gl_j\}_{j=J}^\infty$ be a sequence of nonzero real numbers such that $\lim_{j\to +\infty} \gl_j=0$. Let $\fPhi$ in \eqref{ns:Phi} and $\fPsi^{j}$ in \eqref{ns:Psi} be subsets of distributions in $\Dpr$ for all integers $j\ge J$.
Then the following statements are equivalent
\begin{enumerate}
\item $\FWS_J(\fPhi; \{\fPsi^j\}_{j=J}^\infty)$ is a frequency-based nonstationary tight wavelet frame in $\Lp{2}$, that is, $\fPhi, \fPsi^j\subseteq \Lp{2}$ for all $j\ge J$ and
\begin{equation}\label{ntwf}
\sum_{\ell=1}^\mphi \sum_{k\in \Z} |\la \ff, \fphi^\ell_{\gl_J; 0, k}\ra|^2+\sum_{j=J}^\infty \sum_{\ell=1}^{\mpsi_j} \sum_{k\in \Z} |\la \ff, \fpsi^{j,\ell}_{\gl_j;0, k}\ra|^2=2\pi \|\ff\|^2_{\Lp{2}}, \qquad \forall\; \ff\in \Lp{2};
\end{equation}
\item the pair $(\FWS_J(\fPhi; \{\fPsi^j\}_{j=J}^\infty), \FWS_J(\fPhi; \{\fPsi^j\}_{j=J}^\infty))$ forms a pair of frequency-based nonstationary dual wavelet frames in the distribution space $\Dpr$;
\item $\fPhi, \fPsi^j\subseteq \lLp{2}$ and \eqref{I:ndwf:eq1}, \eqref{I:ndwf:eq2} hold with $\tilde \fPhi:=\fPhi$ and $\tilde \fPsi^j:=\fPsi^j$ for all $j\ge J$.
\item there exist $\phi^1, \ldots, \phi^\mphi, \psi^{j,1}, \ldots, \psi^{j,\mpsi_j}\in \Lp{2}$ for all $j\ge J$ such that $\fphi^1=\wh{\phi^1}, \ldots, \fphi^\mphi=\wh{\phi^\mphi}, \fpsi^{j,1}=\wh{\psi^{j,1}}, \ldots, \fpsi^{j,\mpsi_j}=\wh{\psi^{j,\mpsi_j}}$ for all $j\ge J$, and
\begin{equation}\label{time:ntwf}
\sum_{\ell=1}^\mphi \sum_{k\in \Z} |\la f, \phi^\ell_{\gl_J^{-1}; k}\ra|^2+\sum_{j=J}^\infty \sum_{\ell=1}^{\mpsi_j} \sum_{k\in \Z} |\la f, \psi^{j,\ell}_{\gl_j^{-1};k}\ra|^2=\|f\|^2_{\Lp{2}}, \qquad \forall\; f\in \Lp{2}.
\end{equation}
\end{enumerate}
\end{cor}

For frequency-based nonstationary wavelets with multiresolution-like structure, we have

\begin{prop}\label{prop:ndwf}
Let $J_0$ be an integer and $\{\gl_j\}_{j=J_0}^\infty$ be a sequence of nonzero real numbers such that $\lim_{j\to +\infty} \gl_j=0$. For integers $j\ge J_0$, let $\fPsi^j, \tilde \fPsi^j$ in \eqref{ns:Psi} and
\begin{equation}\label{ns:Phi:2}
\fPhi^j:=\{\fphi^{j,1}, \ldots, \fphi^{j,\mphi_j}\}, \quad
\tilde \fPhi^j:=\{\tilde \fphi^{j,1}, \ldots, \tilde \fphi^{j,\mphi_j}\}
\end{equation}
be subsets of $\lLp{2}$. Then
\begin{equation}\label{ndwf:seq}
(\FWS_J(\fPhi^{J}; \{\fPsi^{j}\}_{j=J}^{\infty}),
\FWS_J(\tilde \fPhi^{J}; \{\tilde \fPsi^{j}\}_{j=J}^{\infty}))
\end{equation}
forms a pair of frequency-based nonstationary dual wavelet frames in the distribution space $\Dpr$ for every integer $J\ge J_0$, if and only if,
\begin{equation}\label{mra:ndwf:eq1}
I^{\gl_j k}_{\fPhi^j}(\gl_j\xi)+I^{\gl_j k}_{\fPsi^j}(\gl_j\xi)=I_{\fPhi^{j+1}}^{\gl_{j+1} k}(\gl_{j+1}\xi), \qquad a.e.\; \xi\in \R, k\in [\gl_j^{-1} \Z]\cup [\gl_{j+1}^{-1}\Z],\;
j \ge J_0
\end{equation}
and
\begin{equation}\label{mra:ndwf:eq2}
\lim_{j \to +\infty} \sum_{\ell=1}^{\mphi_j} \ol{\fphi^{j, \ell}(\gl_j \cdot)}\tilde \fphi^{j,\ell}(\gl_j\cdot)=1 \qquad \hbox{in the sense of distributions},
\end{equation}
where $I^k_{\fPsi^j}$, $k\in \R$, are defined in \eqref{I:Psi} and
\begin{equation}\label{I:Phi}
I_{\fPhi^j}^k(\xi):=\sum_{\ell=1}^{\mphi_j} \ol{\fphi^{j,\ell}(\xi)}\tilde \fphi^{j,\ell}(\xi+2\pi k), \qquad k\in \Z \quad \hbox{and}\quad
I_{\fPhi^j}^k(\xi)\equiv 0, \qquad k\in \R \bs \Z.
\end{equation}
\end{prop}

\begin{proof} By the same argument as in Theorem~\ref{thm:main:1:general}, we see that the pair in \eqref{ndwf:seq} forms a pair of frequency-based nonstationary dual wavelet frames in $\Dpr$ for all integers $J\ge J_0$, if and only if,
\begin{equation}\label{s4:eq4}
\begin{split}
\sum_{\ell=1}^{\mphi_j} \sum_{k\in \Z} &\la \ff, \fphi^{j,\ell}_{\gl_j;0,k}\ra \la \tilde \fphi^{j,\ell}_{\gl_j;0, k}, \fg\ra +
\sum_{\ell=1}^{\mpsi_j} \sum_{k\in \Z} \la \ff, \fpsi^{j,\ell}_{\gl_j;0,k} \ra \la \tilde \fpsi^{j,\ell}_{\gl_j; 0, k}, \fg\ra\\
&=\sum_{\ell=1}^{\mphi_{j+1}} \sum_{k\in \Z} \la \ff, \fphi^{j+1,\ell}_{\gl_{j+1}; 0, k}\ra \la \tilde \fphi^{j+1,\ell}_{\gl_{j+1};0,k}, \fg\ra, \qquad\; \ff, \fg\in \D, j\ge J_0
\end{split}
\end{equation}
and
\begin{equation}\label{s4:eq5}
\lim_{j\to +\infty} \sum_{\ell=1}^{\mphi_j} \sum_{k\in \Z} \la \ff, \fphi^{j,\ell}_{\gl_j;0,k}\ra \la \tilde \fphi^{j,\ell}_{\gl_j;0,k}, \fg\ra=2\pi \la \ff, \fg\ra, \qquad \ff, \fg \in \D.
\end{equation}
By Lemma~\ref{lem:twolevel}, \eqref{mra:ndwf:eq1} is equivalent to \eqref{s4:eq4}. By Lemma~\ref{lem:to1}, \eqref{mra:ndwf:eq2} is equivalent to \eqref{s4:eq5}.
\end{proof}

For $\fa(\xi)=\sum_{k=k_1}^{k_2} a(k) e^{-i k\xi}$ with $a(k_1) a(k_2)\ne 0$, the degree of $\fa$ is defined to be $\deg(\fa):=\max(|k_1|, |k_2|)$.
We finish this paper by the following result which connects a nonstationary wavelet filter bank obtained via a generalized nonstationary oblique extension principle with a pair of frequency-based nonstationary dual wavelet frames in the distribution space.

\begin{theorem}\label{thm:mra:ndwf}
Let $\{\df_j\}_{j=1}^\infty$ be a sequence of nonzero integers such that $\lim_{j\to +\infty} \prod_{n=1}^j |\df_n|=\infty$. Define $\gl_0:=1$ and $\gl_j:=\Big(\prod_{n=1}^j \df_n\Big)^{-1}$ for all $j \in \N$. Let $\fa^j$ and $\tilde \fa^j$, $j\in \N$, be $2\pi$-periodic trigonometric polynomials such that
$\fa^j(0)=\tilde \fa^j(0)=1$ for all $j\in \N$
and
\begin{equation}\label{mask:cond:2}
C:=\sum_{j=1}^\infty |\gl_j| \deg(\fa^j)\|\fa^j\|_{\TLp{\infty}}<\infty \quad \hbox{and}\quad
\tilde C:=\sum_{j=1}^\infty |\gl_j| \deg(\tilde \fa^j)\|\tilde \fa^j\|_{\TLp{\infty}}<\infty.
\end{equation}
Define
\begin{equation}\label{Phi:ndwf}
\fphi^{j}(\xi):=\prod_{n=1}^\infty \fa^{j+n}(\gl_{j+n} \gl_j^{-1} \xi) \quad \hbox{and}\quad
\tilde \fphi^{j}(\xi):=\prod_{n=1}^\infty \tilde \fa^{j+n}(\gl_{j+n} \gl_j^{-1} \xi),
\qquad \xi\in \R, \; j\in \N\cup\{0\}.
\end{equation}
Then all $\fphi^{j-1}$, $\tilde \fphi^{j-1}, j\in \N$ are elements in $\lLp{\infty}$ satisfying
\begin{equation}\label{ns:refeq}
\fphi^{j-1}(\df_j\xi)=\fa^j(\xi)\fphi^j(\xi)\quad \mbox{and} \quad \tilde \fphi^{j-1}(\df_j\xi)=\tilde \fa^j(\xi)\tilde\fphi^j(\xi), \qquad
\forall\; \xi\in \R, j\in \N.
\end{equation}
Let $\fth^{j,1}, \ldots, \fth^{j,\mphi_{j-1}}$, $\fb^{j,1}, \ldots, \fb^{j,\mpsi_{j-1}}$ and
$\tilde \fth^{j,1}, \ldots, \tilde \fth^{j,\mphi_{j-1}}$, $\tilde \fb^{j,1}, \ldots, \tilde \fb^{j,\mpsi_{j-1}}$
with $\mphi_{j-1}, \mpsi_{j-1}\in \N$ and $j\in \N$ be $2\pi$-periodic measurable functions in $\lLp{2}$. Define
\begin{align}
&\fphi^{j-1,\ell}(\xi):=\fth^{j,\ell}(\xi)\fphi^{j-1}(\xi) \quad \mbox{and}\quad \tilde \fphi^{j-1,\ell}(\xi):=\tilde \fth^{j,\ell}(\xi)\tilde \fphi^{j-1}(\xi),\quad j\in \N, \ell=1, \ldots, \mphi_{j-1}, \label{Phi:ns}\\
&\fpsi^{j-1, \ell}(\df_j \xi):=\fb^{j,\ell}(\xi)\fphi^j(\xi) \quad \hbox{and}\quad \tilde \fpsi^{j-1, \ell}(\df_j \xi):=\tilde \fb^{j,\ell}(\xi)\tilde \fphi^j(\xi), \quad j\in \N, \ell=1, \ldots, \mpsi_{j-1}.\label{Psi:ns}
\end{align}
Then $\fPhi^j, \tilde \fPhi^j$ in \eqref{ns:Phi:2} and $\fPsi^j, \tilde \fPsi^j$ in \eqref{ns:Psi} are subsets of $\lLp{2}$ for all $j\ge 0$. Moreover, the pair in \eqref{ndwf:seq}
%
%
forms a pair of frequency-based nonstationary dual wavelet frames in $\Dpr$ for every integer $J\ge 0$, if and only if, for all $j\in \N$,
\begin{align}
&\Theta^j(\df_j\xi)\ol{\fa^j(\xi)}\tilde \fa^j(\xi)+\sum_{\ell=1}^{\mpsi_{j-1}}
\ol{\fb^{j,\ell}(\xi)}\tilde \fb^{j,\ell}(\xi)=\Theta^{j+1}(\xi), \qquad a.e.\, \xi\in \R,\label{mra:filter:eq1}\\
&\Theta^j(\df_j\xi)\ol{\fa^j(\xi)}\tilde \fa^j(\xi+\tfrac{2\pi \omega}{\df_j})+\sum_{\ell=1}^{\mpsi_{j-1}}
\ol{\fb^{j,\ell}(\xi)}\tilde \fb^{j,\ell}(\xi+\tfrac{2\pi \omega}{\df_j})=0, \;\; a.e.\; \xi\in \R,\; \omega\in \Z \bs [\df_j \Z],\label{mra:filter:eq2}
\end{align}
and
\begin{equation}\label{Theta:to1}
\lim_{j\to +\infty} \Theta^{j+1}(\gl_j\cdot)=1\qquad \mbox{in the sense of distributions},
\end{equation}
where
\begin{equation}\label{def:Theta}
\Theta^{j}(\xi):=\sum_{\ell=1}^{\mphi_{j-1}} \ol{\fth^{j,\ell}(\xi)}\tilde \fth^{j,\ell}(\xi), \qquad j\in \N.
\end{equation}
\end{theorem}

\begin{proof} We first establish an inequality on the decay of $\fphi^j$, which plays a critical role in this proof to show \eqref{mra:ndwf:eq2}. By $\fa^j(0)=\tilde \fa^j(0)=1$ and Bernstein inequality $\| [\fa^j]'\|_{\TLp{\infty}} \le \deg(\fa^j) \| [\fa^j]\|_{\TLp{\infty}}$,
\begin{equation}\label{s4:eq6}
|1-\fa^j(\xi)|=\Big| \int_0^\xi [\fa^j]'(\zeta) d\zeta\Big|\le \| [\fa^j]'\|_{\TLp{\infty}} |\xi|\le  |\xi| \deg(\fa^j) \|\fa^j\|_{\TLp{\infty}},
\end{equation}
where $[\fa^j]'$ denotes the derivative of $\fa^j$.  Observe a simple inequality
\begin{equation}\label{ineq}
|z|\le e^{|1-z|}\qquad \forall\; z\in \C.
\end{equation}
In fact, denote $z=re^{i\theta}$ for $r\ge 0$ and $\theta\in \R$. It is easy to prove that $e^{2r}\ge (r+1)^2$ for $r\ge 0$. Consequently, $e^{2|z|}\ge (r+1)^2\ge r^2-2r\cos\theta+1=|z-1|^2$, which leads to \eqref{ineq}.

Now for any $m_1, m_2\in \N\cup\{+\infty\}$ with $m_1\le m_2$, by \eqref{s4:eq6} and \eqref{ineq}, we deduce that
\[
\Big| \prod_{n=m_1}^{m_2} \fa^{j+n}(\gl_{j+n}\xi)\Big|
\le e^{\sum_{n=m_1}^{m_2} |1-\fa^{j+n}(\gl_{j+n}\xi)|}
\le e^{|\xi| \sum_{n=m_1}^{m_2} |\gl_{j+n}| \deg(\fa^{j+n}) \|\fa^{j+n}\|_{\TLp{\infty}}}.
\]
%
Hence, by \eqref{mask:cond:2}, 
we have
\begin{equation}\label{s4:eq7}
\Big| \prod_{n=m_1}^{m_2} \fa^{j+n}(\gl_{j+n}\xi)\Big|\le e^{C|\xi|} \qquad \forall\; \xi\in \R, m_1, m_2\in \N\cup\{+\infty\}\quad \hbox{with}\quad m_1\le m_2.
\end{equation}
On the other hand, by a similar idea as in \cite[Page 93]{Han:siam:2008} and \cite[Page 932]{HanShen:siam:2008}, we have
\[
1-\prod_{n=m_1}^{m_2} \fa^{j+n}(\gl_{j+n}\xi)=\sum_{m=m_1}^{m_2}
\Big(1-\fa^{j+m}(\gl_{j+m}\xi)\Big) \Big( \prod_{n=m+1}^{m_2} \fa^{j+n}(\gl_{j+n}\xi)\Big),
\]
where $\prod_{n=m_2+1}^{m_2}:=1$.
By \eqref{s4:eq6} and \eqref{s4:eq7}, we deduce from the above identity that
\[
\Big |1-\prod_{n=m_1}^{m_2} \fa^{j+n}(\gl_{j+n}\xi)\Big |\le
e^{C |\xi|} \sum_{m=m_1}^{m_2} |1-\fa^{j+m}(\gl_{j+m}\xi)|
\le e^{C |\xi|} |\xi| \sum_{m=m_1}^{m_2} |\gl_{j+m}| \deg(\fa^{j+m})\|\fa^{j+m}\|_{\TLp{\infty}}.
\]
That is, for all $m_1, m_2\in \N\cup\{+\infty\}$ with $m_1\le m_2$, we have
\begin{equation}\label{s4:proda}
\Big |1-\prod_{n=m_1}^{m_2} \fa^{j+n}(\gl_{j+n}\xi)\Big |
\le e^{C |\xi|} |\xi| \sum_{m=j+m_1}^{j+m_2} |\gl_{m}| \deg(\fa^{m})\|\fa^{m}\|_{\TLp{\infty}}\qquad \forall\; \xi\in \R.
\end{equation}
The above inequality implies the uniform convergence of $\prod_{n=1}^\infty \fa^{j+n}(\gl_{j+n}\xi)$ for $\xi$ on any bounded set. Since $\fphi^{j}(\gl_j\xi)=\prod_{n=1}^\infty \fa^{j+n}(\gl_{j+n}\xi)$, we conclude from \eqref{s4:eq7} that
$\fphi^{j}\in \lLp{\infty}$.
Since $\fth^{j,\ell}, \fb^{j,\ell}\in \lLp{2}$,
it is evident that all $\fphi^{j,\ell}, \fpsi^{j,\ell}$ are elements in $\lLp{2}$. Similarly, we can prove
$\tilde \fphi^{j}\in \lLp{\infty}$ and all $\tilde \fphi^{j,\ell}, \tilde \fpsi^{j,\ell}$ are elements in $\lLp{2}$.

Note that $\fphi^{j-1}(\gl_{j-1}\xi)=\fa^j(\gl_j\xi)\fphi^j(\gl_j \xi)$, that is, $\fphi^{j-1}(\xi)=\fa^j(\df_j^{-1}\xi)\fphi^j(\df_j^{-1}\xi)$. So, \eqref{ns:refeq} holds. Also note that $\lim_{j\to \infty} \gl_j=0$ and $\fphi^{j-1}(\xi)=0$ for at most countably many $\xi\in \R$.
Note that \eqref{mra:ndwf:eq1} with $j$ being replaced by $j-1$ is equivalent to
\begin{equation}\label{mra:ns:filter}
\begin{split}
\sum_{\ell=1}^{\mphi_{j-1}} \ol{\fphi^{j-1,\ell}(\df_j\xi)}\tilde \fphi^{j-1,\ell}(\df_j \xi+2\pi \df_j k+2\pi \omega)&+\sum_{\ell=1}^{\mpsi_{j-1}} \ol{\fpsi^{j-1,\ell}(\df_j\xi)}
\tilde \fpsi^{j-1,\ell}(\df_j\xi+2\pi \df_j k+2\pi \omega)\\
&=\delta(\omega)\sum_{\ell=1}^{\mphi_j} \ol{\fphi^{j,\ell}(\xi)}\tilde \fphi^{j,\ell}(\xi+2\pi k+\tfrac{2\pi\omega}{\df_j})
\end{split}
\end{equation}
for all $k\in \Z$ and $\omega=0, \ldots, |\df_j|-1$, where $\delta(0)=1$ and $\delta(\omega)=0$ for all $\omega \ne 0$.
Now it is easy to directly verify that \eqref{mra:filter:eq1} and \eqref{mra:filter:eq2} are equivalent to \eqref{mra:ndwf:eq1}.

We now show that \eqref{mra:ndwf:eq2} is equivalent to \eqref{Theta:to1}.
By $\fphi^{j}(\gl_j\xi)=\prod_{n=1}^\infty \fa^{j+n}(\gl_{j+n}\xi)$ and \eqref{s4:proda}, we have
\begin{equation}\label{s4:eq8}
|1-\fphi^j(\gl_j\xi)|\le e^{C |\xi|} |\xi| \sum_{m=j+1}^\infty
|\gl_{m}| \deg(\fa^{m})\|\fa^{m}\|_{\TLp{\infty}}, \qquad \xi\in \R, j\in \N\cup\{0\}.
\end{equation}
By a similar argument, we also have
\begin{equation}\label{s4:eq9}
|1-\tilde \fphi^j(\gl_j\xi)|\le e^{\tilde C |\xi|} |\xi| \sum_{m=j+1}^\infty
|\gl_{m}| \deg(\tilde \fa^{m})\|\tilde \fa^{m}\|_{\TLp{\infty}}, \qquad \xi\in \R, j\in \N\cup\{0\}.
\end{equation}
By \eqref{mask:cond:2}, we infer that
\[
\lim_{j\to+\infty}
\sum_{m=j+1}^\infty
|\gl_{m}| \deg(\fa^{m})\|\fa^{m}\|_{\TLp{\infty}}=0
\quad\hbox{and}\quad
\lim_{j\to+\infty}
\sum_{m=j+1}^\infty
|\gl_{m}| \deg(\tilde \fa^{m})\|\tilde \fa^{m}\|_{\TLp{\infty}}=0.
\]
By \eqref{Phi:ns}, we deduce that
\[
\sum_{\ell=1}^{\mphi_j} \ol{\fphi^{j,\ell}(\gl_j\xi)} \tilde \fphi^{j,\ell}(\gl_j\xi)=\Theta^{j+1}(\gl_j\xi) \ol{\fphi^j(\gl_j\xi)}\tilde \fphi^j(\gl_j\xi).
\]
Now by \eqref{s4:eq8} and \eqref{s4:eq9}, using Lebesgue dominated convergence theorem, by the same argument as in Theorem~\ref{thm:main:2},
we conclude that \eqref{mra:ndwf:eq2} holds if and only if
\eqref{Theta:to1} holds. By Proposition~\ref{prop:ndwf}, the proof is completed.
\end{proof}


\begin{thebibliography}{10}

\bibitem{BGN}
L. Borup, R. Gribonval, and M. Nielsen,
Bi-framelet systems with few vanishing moments characterize Besov spaces, {\em Appl. Comput. Harmon. Anal.} {\bf 17} (2004),
3--28.

\bibitem{Chui:book}
C.~K.~Chui, An introduction to wavelets. Academic Press, Inc., Boston, MA, 1992.


\bibitem{CHS} C. K. Chui,  W. He and  J. St\"ockler,
Compactly supported tight and sibling frames with maximum vanishing
moments, {\em Appl. Comput. Harmon. Anal.} {\bf 13} (2002),
224--262.

\bibitem{CHS:ntwf}
C.~K.~Chui, W.~He, and J.~St\"ockler, Nonstationary tight wavelet frames. II. Unbounded intervals. \emph{Appl. Comput. Harmon. Anal.}  \textbf{18}  (2005), 25--66.

\bibitem{ChuiShi}
C.~K.~Chui and X.~Shi, Orthonormal wavelets and tight frames with arbitrary real dilations.  \emph{Appl. Comput. Harmon. Anal.} \textbf{9}  (2000),  243--264.

\bibitem{CDF}
A.~Cohen, I.~Daubechies, and J.-C.~Feauveau, Biorthogonal bases of compactly supported wavelets. \emph{Comm. Pure Appl. Math.}  \textbf{45} (1992), 485--560.

\bibitem{CohenDyn}
A.~Cohen and N.~Dyn, Nonstationary subdivision schemes and multiresolution analysis. \emph{SIAM J. Math. Anal.}  \textbf{27}  (1996),  1745--1769.

\bibitem{Daub}
I.~Daubechies, Orthonormal bases of compactly supported wavelets. \emph{Comm. Pure Appl. Math.}  \textbf{41} (1988),  909--996.

\bibitem{Daub:book}
I.~Daubechies, Ten lectures on wavelets, SIAM, CBMS Series, 1992.

\bibitem{DGM}
I.~Daubechies, A.~Grossmann, and Y.~Meyer, Painless nonorthogonal expansions.  \emph{J. Math. Phys.}  \textbf{27}  (1986), 1271--1283.

\bibitem{DH:dwf}
I.~Daubechies and B.~Han, Pairs of dual wavelet frames from any two refinable functions, \emph{Constr. Approx.}, \textbf{20} (2004), 325--352.

\bibitem{DHRS}
I.~Daubechies, B. ~Han, A.~Ron, and Z.~Shen, Framelets: MRA-based
constructions of wavelet frames, {\em Appl. Comput. Harmon. Anal.}
{\bf 14} (2003), 1--46.

\bibitem{E1}
M. Ehler, On multivariate compactly supported bi-frames, {\em J. Fourier Anal. Appl.} \textbf{13} (2007), 511--532.

\bibitem{E2}
M. Ehler, Compactly supported multivariate pairs of dual wavelet frames obtained by convolution, {\em Int. J. Wavelets, Multiresolut. Inf. process.} \textbf{6} (2008), 183--208.

\bibitem{E3}
M. Ehler, Nonlinear approximation associated with nonseparable wavelet bi-frames, {\em J. Approx. Theory}, doi: 10.1016/j.jat.2008.09.007, in press.

\bibitem{FGWW}
M. Frazier, G. Garrig\'os, K. Wang, and G. Weiss,
A characterization of functions that generate wavelet and related expansion, {\em J. Fourier Anal. Appl.} \textbf{3} (1997), 883--906.

\bibitem{Han:msc}
B. Han, Wavelets, M.Sc. thesis at Institute of Mathematics, the Chinese Academy of Sciences, June 1994.


\bibitem{Han:frame}
B.~Han, On dual wavelet tight frames, {\em Appl. Comput. Harmon.
Anal.}, \textbf{4} (1997), 380--413.

\bibitem{Han:jcam}
B.~Han, Compactly supported tight wavelet frames and orthonormal wavelets of exponential decay with a general dilation matrix, \emph{J. Comput. Appl. Math.}, \textbf{155} (2003), 43--67.

\bibitem{Han:sm}
B. Han, Computing the smoothness exponent of a symmetric multivariate refinable function, \emph{SIAM J. Matrix Anal. Appl.}, \textbf{24} (2003), 693--714.

\bibitem{Han:acm:2006}
B.~Han, Solutions in Sobolev spaces of vector refinement equations with a general dilation matrix, \emph{Adv. Comput. Math.}, \textbf{24} (2006), 375--403

\bibitem{Han:siam:2008}
B.~Han, Refinable functions and cascade algorithms in weighted spaces with H\"older continuous masks.  \emph{SIAM J. Math. Anal.} \textbf{40} (2008),  70--102.

\bibitem{Han:oep}
B.~Han, Dual multiwavelet frames with high balancing order and
compact fast frame transform, {\em Appl. Comput. Harmon. Anal.}
\textbf{26}  (2009), 14--42.


\bibitem{HanMo:siamaa}
B. Han and Q. Mo, Splitting a matrix of Laurent polynomials with symmetry and its application to symmetric framelet filter banks, \emph{SIAM J. Matrix Anal. Appl.}, \textbf{26} (2004), 97--124.

\bibitem{HanMo:acha}
B. Han and Q. Mo, Symmetric MRA tight wavelet frames with three generators and high vanishing moments, \emph{Appl. Comput. Harmon. Anal.}, \textbf{18} (2005), 67--93,


\bibitem{HanShen:siam:2008}
B.~Han and Z.~Shen, Compactly supported symmetric $C^\infty$ wavelets with spectral approximation order, \emph{SIAM J. Math. Anal.}, \textbf{40} (2008), 905--938.

\bibitem{HanShen:ca:2009}
B.~Han and Z.~Shen, Dual wavelet frames and Riesz bases in Sobolev spaces, \emph{Constr. Approx.}, \textbf{29} (2009), 369--406.

\bibitem{HW}
E.~Hern\'andez and G.~Weiss, A first course on wavelets. CRC Press, Boca Raton, 1996.


\bibitem{Mallat:book}
S.~Mallat, A wavelet tour of signal processing. Third edition. Elsevier/Academic Press, Amsterdam, 2009.

\bibitem{Meyer:book}
Y.~Meyer, Wavelets and operators. Cambridge University Press, Cambridge, 1992.

\bibitem{RonShen:dwf}
A.~Ron and Z.~Shen, Affine systems in $L_2(\R^d)$ II. Dual systems. {\em J. Fourier Anal. Appl.} {\bf 3} (1997), 617--637.

\bibitem{RonShen:twf}
A.~Ron and Z.~Shen, Affine systems in $L_2(\R^d)$: the analysis of the analysis operator. {\em J. Funct. Anal.}  {\bf 148}  (1997),  408--447.

\end{thebibliography}
\end{document}